\documentclass[11pt,leqno]{amsart}
\makeatletter
\AtBeginDocument{%
  \let\original@@tocwrite=\@tocwrite
  \newif\ifAHVflag
  \def\jlreq@uniqtoken{\jlreq@uniqtoken}
  \def\jlreq@endmark{\jlreq@endmark}
  \long\def\jlreq@getfirsttoken#1#{\jlreq@getfirsttoken@#1\bgroup\jlreq@endmark}
  \long\def\jlreq@getfirsttoken@#1#2\jlreq@endmark#3\jlreq@endmark{#1}
  \renewcommand{\@tocwrite}[2]{%
    \begingroup
      \AHVflagfalse 
      \@ifempty{#2}{}{%
        \expandafter\expandafter\expandafter\ifx\jlreq@getfirsttoken#2\jlreq@uniqtoken{}\jlreq@endmark\Sectionformat\expandafter\@firstoftwo\else\expandafter\@secondoftwo\fi
        {%
          \def\Sectionformat##1##2{\@ifempty{##1}{}{\AHVflagtrue}}%
          #2
        }{\AHVflagtrue}%
      }%
      \def\@tempa{}%
      \ifAHVflag\def\@tempa{\original@@tocwrite{#1}{#2}}\fi
    \expandafter\endgroup
    \@tempa
  }%
}
\makeatother

\usepackage{a4wide}
\usepackage{amssymb}
\usepackage{amsmath}
\usepackage{amsthm}
\usepackage[all]{xy}
\usepackage[usenames,dvipsnames]{xcolor}
\usepackage{lmodern}
\usepackage[T1]{fontenc}

\newcommand{\Hom}{\operatorname{Hom}}
 \DeclareMathOperator{\RHom}{R\underline{Hom}}
\newcommand{\Ext}{\operatorname{Ext}}
\newcommand{\Ind}{\operatorname{Ind}}

\newcommand{\ind}{\operatorname{ind}}

\newcommand{\End}{\operatorname{End}}
\newcommand{\Ord}{\operatorname{Ord}}
\newcommand{\Gal}{\operatorname{Gal}}

\newcommand{\Mod}{\operatorname{Mod}}

\newcommand{\Spec}{\operatorname{Spec}}

\newcommand{\St}{\operatorname{St}}

 \DeclareMathOperator{\Ad}{Ad}

\DeclareMathOperator{\Aut}{Aut}

   \DeclareMathOperator{\red}{red}

\theoremstyle{plain}

\newtheorem{proposition-definition}[theorem]{Proposition-Definition}

\theoremstyle{definition}

\theoremstyle{remark}

\numberwithin{equation}{section}

\title{Representations of $p$-adic  groups over  coefficient  rings}

\author{ Marie-France Vign\'eras}

\date{\today}

\begin{document} 

\maketitle

 \begin{abstract} Motivated by the Langlands program in representation theory, number theory and geometry, the theory of representations of a reductive $p$-adic group  over  a coefficient ring  different from the field of complex numbers  has been widely developped during the last two decades. This article provides  a survey of   basic results obtained  in  the 21st century.
  \end{abstract}

\setcounter{tocdepth}{2}  
\tableofcontents

\section{Introduction}\label{s:0} The theory of representations of a  $p$-adic group $G$, for instance $GL(n,\mathbb Q_p)$, where $\mathbb Q_p$ is the $p$-adic completion of $\mathbb Q$ is an essential part of the  Langlands program. At the beginning, it was a question of studying representations in a complex vector space. But the development of its links with number theory and geometry has required to study continuous representations in vector spaces defined over other fields than $\mathbb C$. There are many possibilities for such a generalization. It is easy to replace $\mathbb C$ by an algebraic closure $\mathbb Q_\ell^{ac}$ of a local field  $\mathbb Q_\ell$, where $\ell$ is a prime different from $p$. The choice of a field isomorphism  $\mathbb C \simeq \mathbb Q_\ell ^{ac}$ identifies continuous complex representations of $G$  and continuous  $\ell$-adic representations.  A more difficult case is that of $\ell=p$ because the topology of a $p$-adic group and of $\mathbb Q_p$ are the same.
  One even considers representations with values not in a vector space, but in a module over some commutative ring   like $\mathbb Z[1/p]$ or $\mathbb Z/p^i \mathbb Z, i\geq 1$. The representations over these different categories of coefficient rings are now essential in the theory of automorphic forms. Their theory  has been widely developped since the beginning of the $21$st century and different versions of the local Langlands correspondance have emerged.  
  
We review  the main basic results  on   representations  over coefficient rings \footnote{that we are aware of, without geometry or derived functors} different from $\mathbb C$. 
  In an attempt to make this paper accessible to readers with a wide  range of backgrounds, we give fairly complete definitions of all terminology. Proofs are omitted,  yet we  give some short indication of the key points, we  cite  sources and we provide examples.  For the theory before 2002,   the reader may consult our  book\footnote{Repr\'esentations $\ell$-modulaires d'un groupe r\'eductif $p$-adique avec $\ell$ diff\'erent de $p$, Birka\"user 1996} and our article in  the proceedings of the Bejing ICM. The subject has  remained confined in  research articles  since these last two decades  and we hope that  this  survey provides a navigable route to the litterature.

 \section{Notation} We work with a triple $(F,G,R)$ where  $F$ is the basic field, $G$ the reductive $p$-adic group, $R$ the coefficient ring. 
We assume that   $F$ is  a local non-archimedean field   of ring of integers $O_F$, uniformizer $p_F$ and    residue field  $k_F$ of characteristic $p$  with $q$ elements, $G$ is the group $\underline G(F)$ of  $F$-points of a connected reductive $F$-group $\underline G$, endowed with  the   topology  generated by the open pro-$p$-subgroups\footnote{called a connected reductive $p$-adic group, but beware that  some authors use this terminology only when $F$ contains $\mathbb Q_p$} and    $R$ is a commutative ring\footnote{a ring is supposed  to have a unit}. 
  
  An $R$-representation $V$ of $G$ will always  be  {\bf smooth} (continuous for the discrete topology on $R$). It is {\bf admissible} if for all open compact subgroups $K$ of $G$,  the $R$-module $V^K$ of vectors  fixed by  $K$ is finitely generated.  
 
 The absolute Galois group $\Gal_E$ of a field $E$ is  the  group of automorphisms of  an algebraic closure $E^{ac}$  fixing $E$. 
 For a  prime number $r$, $\mathbb F_r$  is  a field with $r$ elements,    $\mathbb Z_r$ is the ring of integers in the field $ \mathbb Q_r$ of $r$-adic numbers,    $ \mathbb Z_r^{ac}$  is the ring of integers of $\mathbb Q_r^{ac}$. We always denote by $\ell$ a prime number different from $p$.
   
 \bigskip The parabolic and parahoric subgroups of $G$ play an essential role in the theory of 
 $R$-representations of $G$.  When we will work with them we   will  need  more  notation. 
 
  The parabolic subgroups appear for the first time at the section on parabolic induction. We fix a maximal split torus $T$ of $G$  of $G$-centralizer $Z$ and a minimal parabolic subgroup  $B=ZU$ of  unipotent radical $U$ and opposite $B^{op}=ZU^{op}$. We denote   
  $W_G$ the Weyl group quotient of the $G$-normalizer of $T$ by $Z$, $Z^+\subset Z$ the submonoid of elements contracting $U$ by conjugation, $Z^-$ those contracting $U^{op}$,
  $T^+=T\cap Z^+, T^-=T\cap Z^-$.  The group $G$ is split if $T=Z$ and quasi-split if $Z$ is a torus\footnote{When  $G$ is not quasi-split, $Z$ is not commutative}.  A  {\bf standard parabolic subgroup }  of $G$ is a parabolic subgroup containing $B$, that is,  $P=MN=MB$ with unipotent radical $N$ contained in $U$ and Levi subgroup $M $ containing $Z$. 
  The opposite parabolic subgroup $P^{op}=MN^{op}=MB^{op}$ is not standard.
  
The parahoric subgroups appear for the first time at the section on Hecke algebras. We fix 
  a  {\bf special parahoric subgroup} $K$ of $G$, equal to a parahoric subgroup of $G$  fixing a special point $x_0$  of the apartment of $T$ in the adjoint Bruhat-Tits building of $G$, and a {\bf pro-$p$-Iwahori subgroup} $\tilde J$ of $G$, equal to the maximal open normal pro-$p$ subgroup of 
  the  Iwahori subgroup $J$  of $G$  fixing  the alcove of vertex $x_0$ associated to $B$.  The unique parahoric subgroup  of $T$ is the maximal compact subgroup $T^0= T\cap K=T\cap J$ and the quotient $T/T^0$ is   isomorphic to the group $X_*(T)$ of cocharacters of $T$ via $p_F$.   The  unique parahoric subgroup  of 
 the  connected reductive group  $Z$   is $Z^0=Z\cap K = Z\cap J$, and the quotient $Z/Z^0$ is a commutative  finitely generated group (Thomas Haines and Sean Rostami \cite{HR10}).    For  a  standard parabolic subgroup $P=MN$ of $G$,  $M^0=M\cap K$ is a special parahoric subgroup of $M$ and    $\tilde J_M=\tilde J\cap M$ is a pro-$p$ Iwahori subgroup of $M$. Put $N^0=K\cap N$.

The pro-$p$-Iwahori subgroups of $G$ are all $G$-conjugate, but in general
there are only  finitely many $G$-conjugacy classes of special parahoric subgroups of $G$. 

{\it Examples} There are two conjugacy classes of special parahoric subgroups of  $SL(2, F)$. 

The special parahoric subgroups of  $ GL(n,F)$ are conjugate to $GL(n,O_F)$. 

The (pro-$p$) Iwahori subgroups of  $ GL(n,F)$  are conjugate to the inverse image by the quotient map $GL(n,O_F)\to GL(n,k_F)$ of  the upper (strictly) triangular group of $GL(n,k_F)$.

\section{Change of basic field}  The basic field  $F$ is a finite extension of $\mathbb Q_p$ or of $ \mathbb F_p((t))$,  called a {\bf $p$-adic field} in  characteristic  $0$ or a {\bf local function field} in characteristic  $p$.   Many geometric methods demand   $F$ to be a local functiion field.
For example, the proof by  Bao Chau Ngo\footnote{Fields medal in 2020}   of the   fundamental lemma essential in the Langlands theory, which
 asserts an equality between certain linear combinations of integral orbitals over the Lie algebras of $G$ and of  endoscopic groups. 
On the other hand, for $F$ of characteristic $p$, the harmonic analysis is full of traps,  there are inseparable semi-simple elements, there is no exponential map to pass to the Lie algebra and  $G$ does not contain a  co-compact discrete subgroup (except for type $A$), $G$ is not  a $p$-adic Lie group.

\bigskip  But  the basic field $F$ appears only through the residual field  in many  constructions  (endoscopy, buildings, Iwahori Hecke algebras).  This is a key  to  transfer  properties between  $F$ of different characteristics. For instance,   Jean-Loup Waldspurger \cite{Wa06} proved that the fundamental lemma  for $F$ of characteristic $p$ implies the fundamental lemma  for $F$ of characteristic $0$.  
There is another proof with the  general transfer principle of Cluckers and Loeser   in model theory and motivic integration   \cite{CHL11}, \cite{CGH14}).  In the other direction, the fundamental lemma for the automorphic induction for $GL(n,F)$  proved by Guy Henniart and Rebecca Herb for  $F$ of characteristic $0$ was transfered to  $F$ of characteristic $p$ by Henniart and Bertrand Lemaire \cite{HL06} using {\bf close local fields}.   For a positive integer $m$,  two non-archimedean local fields  are {\bf $m$-closed},  if  their rings of integers modulo  the $m$-th power of their respective maximal ideals are isomorphic. The Deligne-Kazhdan philosophy can be loosely stated as:   the representation theory of Galois groups or of reductive  groups over $m$-close local fields are the same ``up to level m''. For instance, Radhika Ganapathy \cite{G15}  proved that  for two $m$-close local fields $F,F'$ and $\underline G$   split,  the category of  complex representations of  $\underline G(F)$ generated by their  invariants  by the $m$-filtration subgroup of an Iwahori subgroup is equivalent to the same category for representations of $\underline G(F')$. For $\underline G$   not split,  she made sense of a natural connected reductive   group $\underline G'$ over $F'$ associated to $\underline G$,  first when $\underline G$ is quasi-split (an $F$-form of a split group) and then  when   $\underline G$ is general (an inner form of a quasi-split group)  (\cite{G19} 3.A and 5.A).  
 
\bigskip  The local field $\mathbb Q_p$ is a completion of $\mathbb Q$ and $\mathbb Q$ is a globalisation of   $\mathbb Q_p$. The local   is simpler than the global. The ring $\mathbb Z_p$ has only one prime ideal, namely $p\mathbb Z_p$, but the ring $\mathbb Z$ has infinitely many prime ideals. The  absolute Galois group $\Gal_{\mathbb Q_p}$  of  $\mathbb Q_p$ is simple compared to  $\Gal_{\mathbb Q}$. In the same vein,  the  local field  $F$ is the  completion of a  (non-unique) global field\footnote{a global field    is a  finite extension of $\mathbb Q$  or of $\mathbb F_p(T)$} $E$ and $E$ is a globalisation of $F$,   the local group $G$  is  a  localisation the group $H$ of  rational points of a   connected reductive group  over a global field,  and $H$ is a globalisation of $G$\footnote{for $F$ of characteristic $p$ Wee-Teck Gan, Luis Lomeli  \cite{GL18}, for $F$ of characteristic $0$, Shahidi (A proof of Langland's Conjecture on Plancherel measures; Complemenrary Series of $\mathfrak{p}$-adic groups, The Annals of Math., Series 2, Vol.132, 2 (1990), 273-330) when $G$ is quasi-split, implying  the general case as in \cite{GL18}}.
An automorphic irreducible  $\mathbb C$-representation $V_A$ of the adelic group $H(A)$ gives by localisation an irreducible   $\mathbb C$-representation $V$ of $G$ and $V_A$  is  a globalisation of  $V$.  The study of automorphic representations  uses the theory of  representations of  reductive  groups over local fields. In the other direction, some   theorems of  representations of    local  groups  are proved  by embedding the local situation into a global one.

The classical local Langlands correspondence introduced by Langlands in 1967-1970  is a generalization of local class field theory from abelian Galois groups to non-abelian Galois groups.  The absolute Galois group $\Gal_{k_F}$ of the finite field  $k_F$ is topologically generated but the Frobenius $Frob(x)=x^q$,  the  subgroup of elements in $\Gal_F$ with image an integral power of $Frob$ in  the natural quotient map $\Gal_F\to \Gal_{k_F}$ is the {\bf Weil group} $W_F$ of $F$ \footnote{
the kernel $I_F$ of the  quotient map  is an extension of $\prod_{\ell\neq p} \mathbb Z_\ell$ by a pro-$p$ group $P_F$}.  The reciprocity map of local class field theory 
$F^*{\tilde \rightarrow}  W_F^{ab}$   identifies  the irreducible $R$-representations of $GL(1,F)$ with the one-dimensional $R$-representations  of $W_F$  when $R$ is an algebraically closed field. Langlands proposed a parametrisation of the irreducible $\mathbb C$-representations of $G$ in terms of $\mathbb C$-representations of $W_F$, with a purely local characterization.  It is a theorem when $G=GL(n,F)$,  generalized  to representations of $GL(n,F)$ over 
 $R=\mathbb F_\ell^{ac}, \ell\neq p $\footnote{proved when $R=\mathbb C$ by G\'erard Laumon,  Michael Rapoport, and Ulrich Stuhler in 1993 if $F\supset \mathbb F_p((t))$, and if $F\supset \mathbb Q_p$ by Michael Harris and Richard Taylor in 2001,   (Guy Henniart gave another proof), and extended by V.  in 2001 to  $R=\mathbb F_\ell^{ac}, \ell\neq p $}.
The first proofs of local class field theory were global. 
 To-day the proofs of the  local Langlands correspondence for $GL(n,F)$   need  global  arguments  except for $n=2$ and $R=\mathbb C$, which has a local proof (Colin Bushnell and Henniart \cite{BH06}). When  $F\supset \mathbb Q_p$ and $R=\mathbb C$, Peter Scholze \cite{Scho13} gave a new purely local characterization.  The   geometrization of  a   (semisimple) Langlands  correspondence  for all $F,G$  and $R=\mathbb Z_\ell$ for almost all $\ell\neq p$,  by Laurent Fargues and Scholze in 2021,  is entirely local.  A local Langlands correspondence for $GL(n,F), F$ of characteristic $0$, over $R=\mathbb F_p^{ac}$   is a very active research area\footnote{there is yet nothing  for $F$ of characteristic $p$ to my knowledge}.

  \section{Change of coefficient ring}  
  
  Many features  of  complex representations of $G$  use  harmonic analysis  only apparently and can generalize to  representations over other  coefficient   rings. For instance,  
  
  a) The theory of discrete series and tempered complex representations has an algebraic and combinatorial flavour\footnote{the asymptotic behaviour of coefficients may be derived from the central exponents of the Jacquet modules}.   It was extended by Dat \cite{D05}    to an algebraically closed field $R$ of characteristic different from $p$ with a non-trivial valuation.
    
  b)  The proof of the  classification of the irreducible complex representations of  an inner form  of $GL(n,F)$ by  Tadic  for $F\supset \mathbb Q_p$ uses harmonic analysis (the simple trace formula).  Alberto Minguez and Vincent S\'echerre \cite{MS13} gave an algebraic proof  for all $F$ and all  algebraic closed fields  $R$ of characteristic different from $p$.
    
   \bigskip    
  A prime $\ell\neq p$ not dividing the order of a torsion element of  $G$ is called {\bf banal} for $G$\footnote{\cite{DHKM20} Lemma 5.22, Corollary 5.23 for other characterizations}. 
A general principle is that the properties  of  complex representations of $G$  described in purely algebraic  terms transfer to  representations of $G$ over   fields  $R$ of characteristic $0$ or  $\ell$ banal. 

{\it Example}  The banal primes for $GL(m,F)$ are those coprime with $q^{i}-1$ for $1\leq i \leq m$.

 \bigskip   The $R$-representations of $G$ form a locally small  abelian Grothendieck category  $\Mod_R(G)$ \cite{V16b}.  For a commutative ring $S$ which is  an $R$-algebra, the $R$-representations of $G$ are related to the $S$-representations of $G$ by   the scalar extension\footnote{called also base change or induction} 
$S\otimes_R - : \Mod_R(G)\to \Mod_S(G)$ 
 and by   the restriction its  right adjoint:  an $S$-representation is considered as an $R$-representation.   One says that an $S$-representation of $G$  in the image of the scalar extension 
 {\bf descends}  to $R$, or {\bf is defined} on $R$.

\bigskip When $R$ is a field, many properties on admissible irreducible $R$-representations of $G$ still assume $R$ to be  algebraically closed although this is  not necessary. A good tool to check if this  is true is the bijection (Henniart-V.\cite{HV19}, \cite{HV22} section 2):  $$V \mapsto BC (V)$$ 

-  from the isomorphism classes of irreducible admissible $R$-representations  $V$ of $G$,  

- to the Galois orbits\footnote{an orbit under the  the group $\Aut_R(R^{ac})$ of $R$-automorphisms of $R^{ac}$} $BC(V)$ of the isomorphism classes of the  irreducible admissible $R^{ac}$-representations of $G$ defined on a finite extension of $R$.

  $BC(V)$   is the set of  isomorphism classes of the   irreducible subquotients $ V^{ac} $ of     $$R^{ac}\otimes_RV \simeq \oplus ^d\oplus _{V^{ac} \in BC(V)} W(V^{ac}),$$ 
where  $d$ is the reduced dimension of  the division $R$-algebra $ \End_{RG} V$ over its center
   $E_V$, the  length of the  $R^{ac}$-representation $R^{ac}\otimes_RV$ of $G$  is  $d[E_V:R]$, the number of element of  $BC(V)$  is  $[E_V^s:R]$ where $E_V^s$ is  the maximal  separable subextension of $E_V/R$, and 
  $W(V^{ac})$ is an indecomposable $R^{ac}$-representation of $G$ of  irreducible subquotients isomorphic to $V^{ac}$ and of length  $[E_V: E_V^s]$.
Any     $V^{ac}\in BC(V)$ is $V$-isotypic as a  $R$-representation of $G$, and    defined on  a maximal subfield   of  $ \End_{RG} V$ (Justin Trias \cite{Trias20}).

 Any irreducible admissible $R^{ac}$-representation of $G$ is absolutely irreducible and   has a central character by the Schur's lemma.  If  the characteristic of  $R$ is different from $p$,  any irreducible $R^{ac}$-representation of $G$ is admissible and  defined on  a finite extension of $R$ \cite{HV22}.

 \bigskip As $G$ is locally a pro-$p$ group, there is no Haar measure on $G$ with values in a commutative ring $R$ where $p$ is not invertible and   $R$-representations of $G$   present new phenomenona. To understand them  is to-day an open question.  
 
{\it Examples} For a field  $R$ of characteristic $p$, any irreducible $R$-representation $V$ of $G$
 with $\dim_R V^K<\infty$    for some open pro-$p$ subgroup $K$ of $G$, is  admissible (Vytautas Paskunas \cite{Pask04},  a simple proof is given  in  (Henniart \cite{H09}).   As $K$ is a pro-$p$ group, $V^K\neq 0$  when $V\neq 0$, like for  a finite group $G$.
 
 Irreducible  implies admissible when  $G=GL(2, \mathbb Q_p)$ (but not   in general). Indeed, one reduces to  $R=\mathbb F_p^{ac}$; in this case irreducible implies  that the centre acts by  a   character  (Laurent Berger \cite{Be12}) hence is admissible by Barthel-Livne  and Breuil \cite{Br03}.

There exists a non-admissible  irreducible   $\mathbb F_p^{ac}$-representation  of $GL(2,F)$ for an unramified extension    $F$  of $\mathbb Q_p$ (Daniel  Le  \cite{Le19}).
 One does not know  if any infinite dimensional irreducible  non-admissible  $\mathbb F_p^{ac}$- representation   of $G$ has a central character, because its dimension is equal to the cardinal of  $\mathbb F_p^{ac}$ and the classical proof with the Schur's lemma  does not work.
    
\bigskip   It happens that a    property of admissible irreducible  representations of $G$ over a field $R$  transfers  to representations of $G$ over any  coefficient field  of the same characteristic.   
     
{\it Examples}  In  characteristic different from $p$, for  the classification of  cuspidal irreducible   $R$-representations of $G$  by compact induction (Henniart-V.\cite{HV22}). 

  In characteristic $p$,  for   the classification of  non-cuspidal \footnote{cuspidal and supersingular will be defined later}  admissible irreducible $R$-representations  of $G$, for    the classification of  non-supersingular  simple modules of the pro-$p$-Iwahori Hecke $R$-algebra  of $G$  (Noriyuki Abe, Henniart, Florian Herzig and V. \cite{AHHV17}, Henniart-V. \cite{HV19}),   for the existence of a supersingular admissible irreducible  $R$-representation of $G$  when $F\supset \mathbb Q_p$ (Herzig, Karol Koziol and V.\cite{HKV20}).

\bigskip  For a prime $r$\footnote{$\ell$ is reserved for the primes different from $p$, think $r=\ell$ or $p$},   an {\bf $r$-adic representation} of $G$ is a representation  of $G$ on a $ \mathbb Q_r^{ac}$-vector space which is continuous for the  $r$-adic topology on the vector space. In this article,   an $R$-representation of $G$ is supposed  always to be smooth.  A $p$-adic representation of $G$ may be not smooth, but an $\ell$-adic representation of $G$ is smooth if $\ell\neq p$.
A $\mathbb Q_r^{ac}$-representation of $G$ is a smooth $r$-adic representation of $G$. The choice of an isomorphism  $$\mathbb C \simeq \mathbb Q_r^{ac}$$ identifies the complex representations of $G$ and the $\mathbb Q_r^{ac}$-representations of $G$. 

 A {\bf mod $r$  representation}  of $G$  is a  $ \mathbb F_r ^{ac}$-representation of $G$.

 An admissible  $ \mathbb Q_r^{ac}$-representation $V$ of $G$ is called   {\bf integral} if $V$ is defined on a finite extension $E/ \mathbb Q_r$ and $V$ contains an  $G$-stable  $\mathbb Z_r^{ac}$-lattice $L$ \footnote{a free $ \mathbb Z_r^{ac}$-submodule of scalar extension $V$ to $ \mathbb Q_r^{ac}$}, admissible as   a $ \mathbb Z_r^{ac}$-representation of $G$ and descending to   $O_E$\footnote{the ring $O_E$ is principal but not  $ \mathbb Z_r^{ac}$. The definition bypasses this  difficulty}.  The mod $r$ representation
 $\red_r(L)= L\otimes_{\mathbb Z_r^{ac}} \mathbb F_r^{ac}$ of $G$ is called 
 the  reduction modulo $r$ of  $L$.
 
 By the  strong Brauer-Nesbitt theorem (V. \cite{V04a}), an $\ell$-adic  representation $V$ of $G$ of  finite length containing    an admissible  $G$-stable $ \mathbb Z_\ell ^{ac}$-lattice $L$ defined on  some $O_E$ as above, the  $\mathbb Z_\ell ^{ac}[G]$-module $L$ is finitely generated, of reduction $\red_\ell (L)$ of  finite length, and the image  of $\red_\ell (L)$ in the Grothendieck group  of finite length $\mathbb F_\ell ^{ac}$-representations  of $G$ does not depend on the choice of $L$; it is called the   {\bf reduction mod $\ell$} of $V$.
    Two finite length integral  $\ell$-adic representations  of $G$  are said to be {\bf congruent modulo $\ell$} when their  reductions modulo $\ell$ are isomorphic.
     This does not hold true  for    $\mathbb Q_p^{ac}$-representations  of $G$.  

{\it Example} An irreducible $\mathbb Q_p^{ac}$-representation $V=\ind_K^GW$ of $G=PGL(n,F)$ compactly induced from  a representation $W$ of $K=PGL(n,O_F)$ contains an admissible  $G $-stable  $\mathbb Z_p^{ac}$-lattice  $L$ defined on  some $O_E$ as above,  of infinite length reduction, and another one $L'$ of finite length reduction.  Take 
 $L= \ind_K^G W_{\mathbb Z_p^{ac}}$ for a
$K$-stable   $\mathbb Z_p^{ac}$-lattice $W_{\mathbb Z_p^{ac}}$ of $W $ and 
$L'=V\cap  \ind_\Gamma^G 1_{\mathbb Z_p^{ac}}$ for   a small enough discrete cocompact subgroup $\Gamma$ of $G$.

\section{Parabolic induction}
For any $F,G,R$ and any parabolic subgroup $P$ of $G$ of Levi quotient $M$, 
the {\bf parabolic induction} \footnote{$\ind_P^G(W)$ is  the right translation of $G$ on the $R$-module of locally constant functions $f:G\to W$  such that $f(mng)=m f(g)$ for $m\in M,n\in N, g\in G$}
$$\ind_P^G:\Mod_R(M)\to \Mod_R(G)$$
 allows to construct $R$-representations of $G$ from $R$-representations of  the smaller connected reductive $p$-adic group $M$. 
  The parabolic induction   has excellent  properties,   it
   commutes with small  direct sums  \cite{V16b} \footnote{when $R$ is a field of characteristic $p$, $ \ind_P^G$ commutes with direct products \cite{SS22}} Lemma 4.3;    for  $p$ nilpotent in $R$, it   is  fully faithful  \cite{V16b}; for a field $R$,  the parabolic induction respects finite length representations with admissible subquotients     (this depends on the classification of admissible irreducible representations  if   the characteristic of $R$  is $p$).

 The parabolic induction is exact and  has   a  {\bf left adjoint}  $L_P^G$ called the Jacquet functor,   equal to the coinvariant functor $ (-)_N$ with respect to the unipotent radical $N$ of $P$,  and  a {\bf right adjoint} \footnote{by  \cite{KS06} 8.3.27 as $\Mod_R(G)$ is a locally small  abelian Grothendieck category anf $\ind_P^G$ is right exact and commutes with small  direct sums} $R_P^G$ \cite{V16b}.  By adjointness, $L_P^G$ is   right exact  and  $R_P^G$ is left exact.  The scalar extension commutes with the three parabolic functors \cite{HV19}.

\bigskip 
 For $p$  invertible in $R$,  the  {\bf second adjunction}   $$R_P^G=  \delta_P L_{P^{op}}^G.$$ where $\delta_P$ is the modulus of $P$\footnote{$ \delta_P (m)= | \det \Ad _{Lie N}(m)|\in q^{\mathbb N}$},
 is a deep property 
 proved this year by  Dat, David Helm, Robert Kurinczuk and   Gilbert Moss (\cite{DHKM22}  Corollary 1.3),  originally proved by  Bernstein  when $R=\mathbb C$.
 When $R$ is noetherian, the parabolic induction $\ind_P^G$ respects admissibility, 
the second adjunction implies  that 
     $\Mod_R(G)$ is   noetherian,  that the parabolic induction respects projective (resp. finitely generated) $R$-representations (\cite{DHKM22}  Corollaries 1.4, 1.5), and that $L_P^G$  respects admissibility.
     The functor $L_P^G$  is exact, preserves infinite direct sums  \cite{D09}, and  when $R$ is a field, $L_P^G$ respects   finite length  because $L_P^G$ respects the property of being finitely generated, and      an admissible finitely generated $R$-representation of $G$  has finite length  (the proof uses  the Moy-Prasad unrefined types when $R$ is algebraically closed but  algebraically closed is not necessary).

 \bigskip  For a field $R$ of characteristic $p$, the adjoint functors 
   $L_P^G$ and $R_P^G$ send  an  admissible  irreducible $R$-representation  of $G$ to $0$ or to an admissible irreducible $R$-representation of $M$.  Irreducible is necessary,   an example of an admissible $R$-representation $V$ of $G$ with  $L_P^G(V)$  not admissible is given in  \cite{AHV19}.  But contrary to the case   $R=\mathbb C$, the functors $L_P^G$ and $R_P^G$ fail to be  exact (for $R_P^G$  (Emerton \cite{Em10b}, Koziol \cite{Koz22}), 
  $\ind_P^G$ does not preserve finitely generated representations, $R_P^G$ does not preserve infinite direct sums (Abe-Henniart-V. \cite{AHV19} section 4.5).

\bigskip When  $p$ is nilpotent in the commutative ring $R$, the right adjoint $R_P^G$ respects admissibility (Abe-Henniart-V.  \cite{AHV19}); it  is equal to the Emerton's functor $\Ord_{P^{op}}^G$ of ordinary parts  on admissible $R$-representations\footnote{there is no description of $R_P^G$  on non-admissible  representations}.
If moreover $R$ is artinian,   Matthew Emerton \cite{Em10b} extended the functor of ordinary parts to a $\delta$-functor, expected to coincide with the derived functors  when  the characteristic of $F$ is $0$.

{\it Example} When $G=SL(2,\mathbb Q_p)$,  Koziol \cite{Koz22} showed that the derived functors of $R_B^G$ and $\Ord_B^G$  are equal  on any  absolutely irreducible  $\mathbb F_p^{ac}$-representations of $G$.

 When  the characteristic of $F$ is $p$, suprisingly  $R_P^G$ is exact on admissible  $\mathbb F_p^{ac}$-representations  of $G$ (Julien Hauseux \cite{Hau18b}).

 \bigskip A representation of $G$ over  a  field $R$  is called {\bf unramified} when it is trivial on the subgroup $G^0$ of $G$ generated by its compact subgroups\footnote{This coincides with the classical definition (Henniart-Lemaire \cite{HL17} 2.12 Remarque 1)} .   The group  $\Psi_R(G)$ of  unramified $R$-characters $\psi:G\to R^*$ of $G$  is a torus.   
    {\bf Generic irreducibility} says that
for any parabolic subgroup $P$ of $G$ of Levi $M$ and any irreducible $R$-representation $W$ of $M$, the set of  $\psi\in \Psi_R(M)$ such that $\ind_P^G(W\otimes \psi$) is irreducible is Zariski-dense in  $\Psi_R(M)$.
    Generic irreducibility  is probably true for any $F,G$ and any field $R$.
    
    {\it  Example} Generic irreducibility  is  known  for $R$ of characteristic $p$ (Abe-Henniart-V.\cite{AHV19})  or  when $F\subset \mathbb Q_p$  for  $R$ algebraically closed of characteristic different from $p$ (Dat \cite{D05}).

\bigskip Dat (\cite{D05} Theorem 3.11) extended the {\bf Langlands  quotient theorem} when $R=\mathbb C$ to  any algebraically closed  field $R$ of characteristic different from $p$  with a non-trivial valuation $\nu$ (for example $\mathbb Q_\ell ^{ac}$):

 When $P=MN$ is a standard parabolic subgroup of $G$, $W$ is  a $\nu$-tempered irreducible $R$-representation of $M$,  and $\psi\in \Psi_R(M)$  satisfies  $-\nu(\psi)\in (\mathcal A_P^*)^+$, the $R$-representation $\ind_P^G(W\otimes \psi)$ has a unique irreducible quotient $J(M,W, \psi)$. Any  irreducible $R$-representation $V$ of $G$ is isomorphic to $J(M,W, \psi)$ for a unique  triple 
 $(P,W, \psi)$.

An admissible $R$-representation $V$ of $G$ is {\bf $\nu$-tempered} (Dat \cite{D05} Definition 3.2) if  
  for any standard parabolic subgroup $P=MN$ such that $L_P^G(V)\neq 0$,   any exponent $\chi$ in $L_P^G(V)$ satisfies $-\nu (\delta_P^{-1/2}\chi) \in 
  \overline {{}^+\mathcal A^*_P}$ \footnote{
 Let 
  $\Delta(M)$ denote  the set of simple roots of $T$ in $M$,  $\Delta(P)$
  the set of simple roots in $P$ of  $T_M$, 
   $\mathcal A^*= X\otimes_\mathbb Z \mathbb R$ where $X$ is the lattice of rational characters of   $T$,  $( \ , \ )$ a $W_G$-invariant  scalar product on  $\mathcal A^*$. Then  $\overline {{}^+\mathcal A^*_P} = \sum _{\alpha \in \Delta_P }\mathbb R_{\geq 0}\, \alpha$ and $(\mathcal A^*_P)^+ $ is the cone $\{x  \in \mathcal A^* , \ 
  (x, \alpha)=0  \text{ for } \alpha \in \Delta(M), \ (x, \alpha)>0  \text{ for } \alpha \in \Delta(P)\}$}.  It is called {\bf discrete} if $-\nu (\delta_P^{-1/2}\chi) \in 
 {}^+\mathcal A^*_P$.
   The {\bf exponents} of  $L_P^G(V)$ are the $R$-characters  of the split component $A_M$ of the center of $M$   appearing in $L_P^G(V)$ seen as an $R$-representation of $A_M$.

 From the  Dat's theory  of $\nu$-tempered representations, one deduces  (David Hansen, Tasho Kaletha  and Jared Weinstein (\cite{HKW22} C.2.2)):

The Grothendieck group  of finite length $ \ell$-adic representations of $G$ is generated by representations of the form $\ind_P^G(W  \otimes \psi)$,  for a standard parabolic subgroup $ P\subset G$  of Levi $M$,  an integral irreducible  $ \ell$-adic representation $W$ of $M$  and an unramified $\ell$-adic character $\psi$  of $M$.

\bigskip  
 \section{Admissible representations and duality}  The classification of irreducible admissible $R$-representations of $G$ is  an objective of the local Langlands program. There are few   finite dimensional representations when $G$ is not compact modulo the center $Z(G)$, and admissibility is  a  crucial finiteness property.

When $R$ is a noetherian  commutative ring, a  subrepresentation of  an admissible  $R$-representation of $G$ is admissible. 
A  quotient of an admissible  $R$-representation of $G$ is admissible \cite{V11}  and the category  $\Mod_R(G)^a$ of admissible $R$-representations of $G$ is abelian if  $p$ is  invertible in $R$,  or if  $R$ is a finite field of characteristic $p$ and  $F\supset \mathbb Q_p$\footnote{the completed group algebra of $R[K]$ is noetherian when $F\supset \mathbb Q_p$ but not when $F\supset \mathbb F_p((T))$}.

{\it Example}  When $F\supset \mathbb F_p((T))$  and $p $ is not invertible in $R$,  rhere exists an admissible representation with a non-admissible quotient   
(Abe-Henniart-V.  \cite{AHV19}).

\bigskip   
Let $R$ be a field. The  {\bf smooth dual} $V^\vee$ of an $R$-representation  $V$ of $G$ is the   smooth part of the  contragredient action of $G$ on the   linear dual $V^*=\Hom_R(V,R)$  \footnote{the smooth dual is the set of linear forms on $V$ fixed by some open subgroup of $G$}.  

 For $R$ of  characteristic   different from $p$, the smooth dual is  an autoduality on   $\Mod_R(G)^a$.  In particular, $V^\vee$ is  irreducible if and only if $V$ is irreducible. The smooth dual and the parabolic induction and its left adjoint satisfy \footnote{the normalized induction $\ind_P^G(W\otimes \delta_P^{1/2})$  commutes with the smooth dual, the second isomorphism is equivalent to the second adjunction}:
$$(\ind_P^G W)^\vee \simeq \ind_P^G(W^\vee \delta_P), \  \  L_P^G(V^\vee ) \simeq (L_{P^{op}}^G(V))^\vee ,$$
for any  $R$-representation $W$ of $M$ and any  admissible $R$-representation $V$ of $G$.  


 For $R$ of  characteristic   $p$, 
 the smooth dual of  any infinite dimensional admissible irreducible $R$-representation  of $G$ is zero !
For $F$ of characteristic $0$,  Jan Kohlhaase \cite{Kohl17} developped a  higher smooth duality theory on $\Mod_R(G)^a$.  He  studied  the  $i$-th smooth duality functors $S^i: \Mod_R(G)^a\to \Mod_R(G)^a $ for $0\leq i \leq d $ for $d=\dim_{\mathbb Q_p}G$  under tensor product, inflation and induction and proved that  for $V\in \Mod_R(G)^a$,  the integer $$d(V) = \max \{ i \ | \ S^i(V)\neq 0)\} $$ 
satisfies 

(i)  $d(V)=0$  if and only if $V$ is finite dimensional,  

(ii)    $d(\ind_P^GW)= d(W)+ \dim_{\mathbb Q_p} N$, for   a parabolic subgroup   $P=MN$ and $W\in \Mod_R(M)^a$, 

 (iii) $d(V)=1$ and 
  $S^1(V)$ coincides with  the  Colmez's contragredient  introduced for the $p$-adic Langlands correspondence  for $G=GL(2,Q_p), R=\mathbb F_p^{ac}$ and $V$ irreducible of   infinite dimensional;  for the Steinberg representation $\St_G$  which is irreducible, $S^1(St_G)$ is indecomposable of length $2$ ! 

  For $G$  unramified\footnote{$G$ is quasi-split and  splits over some unramified extension of $F$}, $K$ an hyperspecial subgroup  of $G$,    $W\in \Mod_{\mathbb F_p^{ac}}(K)$ and $i>\dim_{\mathbb Q_p}U$, we have $S^i(\ind_K^GW)=0$ (Claus Sorensen \cite{So19}).  
\section{Supercuspidal    support}
 
  An $R$-representation  $V$ of $G$ is called {\bf cuspidal} if it is killed $$L_P^G(V)=R_P^G(V)=0$$by  the left and right adjoints of the  parabolic induction  for all parabolic subgroups  $P\neq G$.

\bigskip When  $p$ is invertible in $R$, 
the second adjunction implies that $V$ is cuspidal  if and only if  $L_P^G(V)=0$ for any proper parabolic subgroup  $P $ of $G$.  
Any irreducible   $R$-representation $V$ of $G$ is a   subrepresentation of  $\ind_P^GW$ for some cuspidal   irreducible  $R$-representation $W$. Assuming that $R$ is an algebraically closed field\footnote{algebraically closed is probably not necessary},   the pair  $(M,W)$ is  unique modulo  $G$-conjugation, its 
 $G$-conjugation  class  of $(M,W)$ is called the  {\bf cuspidal support} of $V$. Twisting the cuspidal support by unramified characters we get the {\bf inertial cuspidal support} $\Omega$ of $V$. So,    
$\Omega$ is  the set of $(M', W')$ for $(M', W')$ $G$-conjugate to $(M,W\otimes \Psi)$ and $\psi \in \Psi_R(M)$. It  is an algebraic variety with regular functions  
$\mathcal O(\Omega)= (R[M/M^0]^S)^H$ where $S$ is the (finite) group of  $\psi \in \Psi_R(M)$ such that $W\otimes\psi\simeq W$. The  subgroup of   $w\in W_G$ fixing $M$ acts on the $R$-representations of $M$, and  $H$ is the group of those $w$ such that $W^w \simeq W\otimes\psi$ for some $\psi\in \Psi_R(M)$.

When $p$ is not invertible in $R$, one  needs both $L_P^G$ and $R_P^G$ to define cuspidality. The mod $p$ Steinberg representation $\St_G$  and the trivial representation $1_G$ of $G$  satisfy for any parabolic subgroup $P$ of Levi $M$,
$$L_P^G(\St_G)=0, R_P^G(\St_G)=\St_M,\ \ L_P^G(1_G)=1_M, R_P^G(1_G)=0.$$
For a field $R$ of characteristic $p$,  the  Steinberg representation  is not   a   subrepresentation of  $\ind_P^GW$ for any cuspidal   admissible irreducible  $R$-representation $W$. Yet, any irreducible $R$-representation $V$ of $G$ is a subquotient of $\ind_B^GW$ for some $R$-representation $W$ of $Z$ (Abe-Henniart-Herzig-V.\cite{AHHV17} IV.1 for $R$ algebraically closed field).

\bigskip  An  admissible  irreducible $R$-representation   of $G$ which is  not isomorphic to a subquotient of  a proper parabolically induced representation $\ind_P^G W$ for all  $P\neq G, W$ an admissible irreducible $R$-representation of $M$, is  called {\bf supercuspidal}\footnote{One does not need to suppose $W$  irreducible  when $R$ is an algebraically closed field tof characteristic different from $p$ (Dat\cite{D18b})}.

Any admissible irreducible  $R$-representation $V$ of $G$ is a   subquotient  of  $\ind_P^GW$  for some supercuspidal  admissible   irreducible $R$-representation $W$. 

\bigskip  For a field $R$ of  characteristic $p$,   $(P,W)$ is unique modulo $G$-conjugation.  This follows from the  classification.  Any cuspidal  admissible irreducible  $R$-representation  of $G$ is supercuspidal. 
 
\bigskip For a field $R$ of  characteristic different from $p$, a cuspidal   irreducible  $R$-representation $V$ of $G$ is not always supercuspidal.
  The  $G$-conjugation class   of $(M,W)$ is called a {\bf supercuspidal support} of $V$.  When $R$ is algebraically closed, its  twist  by unramified characters is called an {\bf inertial supercuspidal support} of $V$;  if all the irreducible $R$-representations of $G$ have a unique supercuspidal support,  the {\bf Bernstein variety}   $\mathcal B_R(G)$ is the disjoint union of  the inertial  supercuspidal supports of the    irreducible $R$-representations of $G$. 
Contrary to the cuspidal support, 
the supercuspidal support is not always unique. 

{\it Examples} When $G=GL(2,\mathbb Q_p)$,$R=\mathbb F_\ell^{ac}, \ell $ divides $p+1$, the unique  infinite dimensional irreducible subquotient of the representation $\ind_B^G 1_Z$  indecomposable of  length $3$ is cuspidal and 
 non-supercuspidal. 

The supercuspidal support is not always unique when $R=\mathbb F_\ell^{ac}, \ell$ divides $q^2+1$ and $G$ is the finite group $ Sp_8(\mathbb F_q)$ (Olivier Dudas  \cite{Du18}) or $Sp_8(F)$ (Dat  \cite{D18b}).

The supercuspidal support is unique    if  $R$ has characteristic $0$,  or   $G$ is an inner form of $GL(n,F)$  (Minguez-S\'echerre \cite{MS14a}), or  $G$ is the unramified unitary group $U(2,1)$, $p\neq 2$  (Kurinczuk \cite{K14}), when  $R$ is algebraically closed (probably this is not necessary).

\bigskip An irreducible $\mathbb Q_\ell^{ac}$-representation of $G$ is integral if and only if its supercuspidal support is integral (Dat-Helm-Kurinczuk-Moss \cite{DHKM22} Corollary 1.6). Is any irreducible  mod $\ell$ representation  of $G$ a subquotient of the reduction modulo $\ell$ of an integral irreducible $\ell$-adic representation  ?

\bigskip For a field  $R$ of characteristic $0$ or $\ell$ banal for $G$, any  cuspidal irreducible  $R$-representation  of $G$ is supercuspidal, and projective in the category of $R$-representations of $G$ with a given central character.
   The reduction modulo $\ell$  of any integral  cuspidal irreducible $\ell$-adic  representation  of $G$ is irreducible and  cuspidal. The reduction modulo $\ell$ of an integral irreducible $\ell$-adic representation of $G$ may be reducible. Does
  any  irreducible  mod $\ell$ representation  of $G$  lifts to an  irreducible  $\ell$-adic representation of $G$  ? \footnote{is the reduction modulo $\ell$  of an  integral  cuspidal irreducible $\ell$-adic  representation  of $G$}

 \section{Hecke algebras} Hecke $\mathbb Z$-algebras   appear everywhere in   the theory of representations of $G$ to find  algebraic proofs of properties of   representations   proved formerly with harmonic analysis.  
 An open subgroup $K$ of $G$ which is compact  or compact modulo the  center of $G$, defines a {\bf  Hecke ring} 
    $$\mathcal H(G,K) = \End_{\mathbb Z[G]}  \mathbb Z[K\backslash G].$$
naturally isomorphic to the opposite of    $  \mathbb Z[K\backslash G/K]$.
  For any commutative ring $R$, the Hecke $R$-algebra  $\mathcal  H_R(G,K) =\End_{R[G]} R[K\backslash G] $ is the scalar extension to $R$ of the Hecke ring. 
  
  \bigskip   A famous finiteness theorem of Deligne-Bernstein  when $R=\mathbb C$ extended by Dat-Helm-Kurinczuk-Moss \cite{DHKM22}   is  the key of the proof of the second adjunction:

When $R$ is any noetherian $\mathbb Z_\ell$-algebra, the  center  $\mathcal Z_R(G,K) $ of  $\mathcal H_R(G,K)  $ is a finitely generated $R$-algebra and  $\mathcal H_R(G,K)$  is a finitely generated  $\mathcal Z_R(G,K) $-module.

One  proves   an equivalent statement, involving the {\bf Bernstein center} $\mathcal  Z_R(G)$, which is the  endomorphism ring   of the identity functor  of $\Mod_R(G)$: 

 For $R$ as above, any finitely generated $R$-representation $V$ of $G$ is $\mathcal  Z_R(G)$-admissible and the natural image of $\mathcal  Z_R(G)\to \End_{R[G]}V$
is a finitely generated $R$-algebra.

 The highly non-trivial  proof uses the  Fargues-Scholze local version of the Vincent Lafforgue's theory of excursion operators \cite{FS21}.

\bigskip   The $R$-representations of $G$ and  the right $\mathcal H_R(G,K) $-modules are related   by 
  the  {\bf $K$-invariant  functor}:  $$V\mapsto V^K\simeq \Hom_{R[G]}(R[K\backslash G], V): \Mod_R(G) \to \Mod \mathcal H_R(G,K), $$
and by its left adjoint $\mathcal M \to \mathcal M \otimes_{ \mathcal H_R(G,K)}R[K\backslash G]$. From now on, a module of a Hecke algebra will be a right module.  

 When $R$ is a field of  characteristic different from $p$, any simple $ \mathcal H_R(G,K)$-module is finite dimensional  if  the $K$-invariant functor induces a bijection between the 
(isomorphism classes of)  irreducible $R$-representations $V$ of $G$ with $V^K\neq 0$ and the (isomorphism classes of)   $ \mathcal H_R(G,K)$-modules, as irreducible implies admissible.  This is the case  when the order of any finite quotient of $K$ is invertible in $R$ (\cite{HV22} Theorem 3.2), or  when $K$ is an Iwahori group $J$ or a  pro-$p$ Iwahori subgroup $\tilde J$.

When any subrepresentation of any $R$-representation  of $G$ generated by its $K$-invariants has the same property, 
the $K$-invariant  functor gives an equivalence: $$ \Mod_R(G)(K)\tilde{\to} \Mod \mathcal H_R(G,K) $$
By a classical result of Borel,
 the category  of complex representations generated by their vectors  invariant by an Iwahori subgroup $J$ is  an indecomposable factor of $\Mod_{\mathbb C}G)$,  equivalent to $\Mod {\mathcal H}_{\mathbb C}(G,J)$ by the $J$-invariant functor.  The same is true with a pro-$p$ Iwahori subgroup $\tilde J$ without ``indecomposable''. 
 
  There are no equivalences for $R $ of characteristic $p$ and 
 $K=\tilde J$  in general.

{\it Example}  However, the equivalence  is true  for 
  $R=\mathbb F_p^{ac}$ and 
 $K=\tilde J$ if   $G=GL(2, \mathbb Q_p)$ or $SL(2, \mathbb Q_p), p\neq 2$  (Rachel Ollivier \cite{O09}\footnote{supposing that a uniformizer of $F$ acts trivially}, Koziol \cite{Koz16b},  Ollivier-Peter Schneider \cite{OS18}) 

\bigskip For a prime $r$,  a  $\mathbb Q_r$-representation $V$ of $G$  is called {\bf locally integral} if for some finite extension $E/\mathbb Q_r$, $V^K$ admits a $\mathcal H(G,K)$-stable  $O_E$-lattice for all open compact subgroups $K$ of $G$. 

An  integral irreducible  $\mathbb Q_r^{ac}$-representation    is clearly locally integral.  The two notions for irreducible representations coincide if $r=\ell \neq p$  \cite{D05}. 

The  equivalence between integral and locally integral irreducible $\mathbb Q_p^{ac}$-representations of $G$  is an open question.  When $G$ is split, it  is the analogue of the  Breuil-Schneider conjecture \cite{BS07}  restricted to smooth representations (Hu \cite{Hu09}, Sorensen \cite{So13},\cite{So15}).
 A  finite length $\mathbb Q_p^{ac}$-representation  $V$ of $G$ is locally integral if and only if (Dat \cite{D09a}):
 $$\nu( \delta_P^{-1/2} \chi)\in \rho_P - \overline{ {}^+ \mathcal A_P^*}$$
 for any  standard parabolic subgroup $P=MN$ of $G$ with  $L_P^G(V)\neq 0$, and any exponent $ \chi$ of $ L_P^G(V)$\footnote{$\rho_P$ is  half the sum of the roots of $A_M$ in $Lie P$. The formula can  be  simplified ! }.  
  
\section{Representations in characteristic  different from $p$} 
For any commutative ring $R$ and any open subgroup $K$ of $G$,  an $R$-representation $W$ of $K$ define a $R$-representation $\ind_K^GW$ of $G$ by {\bf compact induction}\footnote{ the $R$-module of functions $f:G\to W$  supported on finitely many cosets $Kg$, satisfying $f(kg)=\rho(k)f(g)$ for $k\in K, g\in G$ where $G$ acts by right translation}. 
 
 {\it  Example}  $\ind_K^G1_K=R[K\backslash G]$ for the trivial $R$-representation $1_K$ of $K$.

\bigskip When $R$ is a field of characteristic  different from $p$, 
all cuspidal irreducible $R$-representations of $G$ are conjectured to be compactly induced from open subgroups  of $G$ compact modulo the center of $G$. 

For $R$ algebraically closed, the conjecture has been proved for the level $0$\footnote{definition in  the section on Bernstein blocks}  representations of any $G$  or when

 $G$ has rank $1$ (Martin Weissman \cite{Weiss19}), 
 
 $G$ is an inner form of $GL(n,F)$ (Minguez-S\'echerre \cite{MS14a}), or  of $SL(n,F)$ (Peyi Cui \cite{C19},\cite{C20}),
 
 $G$ is a classical group and $p\neq 2$ (Shaun Stevens \cite{Ste08}, Stevens-Kurinczuk-Daniel Skodlerak \cite{KSS21}) or a quaternionic form of $G$ (Skodlerak \cite{Sko19}), 
 
 $G$ splits on a moderately ramified extension of $F$ and $p$ does not divide the order of the absolute Weyl group  (Jessica Fintzen \cite{F19}),

\bigskip  Algebraically closed is not necessary and there is an explicit 
    list $\mathcal X$ of pairs $(K,W)$ of $G$ where
 $K$ is   an open subgroup  of $G$ compact modulo the center and $W$ an $R$-representation of $K$ such that  $\ind_K^GW$  is irreducible  cuspidal satisfying (Henniart-V.  \cite{HV22}):
 
a) any cuspidal    irreducible $R$-representation of $G$  is isomorphic to $\ind_K^GW$ for some $(K,W)\in \mathcal X$ unique modulo $G$-conjugation,
  
b) $\ind_K^GW$ and  $W$ have the same interwinning algebra $\End_{R[K]}W\tilde{\rightarrow}\End_{R[G]} \ind_K^GW$,
 
c)  $\ind_K^GW$ is supercuspidal if and only if $ W$ is supercuspidal, for the  ``natural notion of   supercuspidality''  of  $W$ \footnote{Fintzen gave another proof when  $G$  is moderately ramified  and $p$ does not divide the order of the absolute Weyl group},

 d) $\mathcal X$ is    stable by  automorphisms of $R$.
 
 \bigskip Until  the end of this section,  $R$  is an algebraically closed field of characteristic  different from $p$ and  
 $G=GL(m,D)$   where $D$ is a central division algebra of dimension $d^2$ over $F$, $n=md$. 
 
 \bigskip  Minguez and S\'echerre    classified the irreducible  $R$-representations  of $G$ with a given supercuspidal support   by ``supercuspidal multisegments'', and those  with a given  cuspidal support  by ``aperiodic cuspidal multisegments''   \cite{MS14b}. This is the  generalisation the    Bernstein-Zelevinski classification of complex irreducible representations of $GL(n,F)$. 
 For $R$ of characteristic $\ell$, the proof uses the theory of $\ell$-modular types  (Minguez-S\'echerre \cite{MS14a}) and deep results on affine Hecke algebras  of type $A$ at roots of unity.

  Any irreducible $\ell$-modular representation of $G$ is a subquotient of the reduction modulo $\ell$ of an integral irreducible $\ell$-adic representation \cite{MS14b}. In the other direction,
 any irreducible $\ell$-modular 
representation $V$ of $G$ lifts to an $\ell$-adic representation  when it is supercuspidal or ``banal'' or unramified\footnote{$V^{GL(m,O_D)}\neq 0$, equivalent to $V$  irreducibly parabolically induced from an unramified character of a Levi subgroup  \cite{MS14c},} (Dat\cite{D05}, Minguez-S\'echerre \cite{MS14b},  \cite{MS13}, \cite{MS14c}) or when $G=GL(n,F)$.
Contrary to the case $G=GL(n,F)$, some irreducible cuspidal $\ell$-modular 
representation of $G$  may not lift and the reduction modulo $\ell$ of a integral cuspidal irreducible $\ell$-adic-representation of $G $ may be reducible; its
 irreducible components are cuspidal  and in the same inertial class.
 
 {\it Example} When  $q=8,\ell=3, d=2$, any integral  irreducible $\ell$-adic representation  of $D^*$ containing an homomorphism $\chi:O_D^*\to (\mathbb Q_\ell^{ac})^*$ trivial on $1+P_D$ such that $\chi\neq \chi^q$ has dimension $2$ and  a non-irreducible  reduction modulo $\ell$. When  $q=4, \ell=17, d=2$,  there exists an irreducible cuspidal $ \ell$-modular representation of  $GL(2,D)$ not lifting to $\mathbb Q_\ell^{ac}$  (Minguez-S\'echerre \cite{MS17}).  

\bigskip Let ${\mathcal D}_{\mathbb C}(G )$ denote the set of isomorphism classes of   the essentially square integrable irreducible (or discrete series) complex representations of $G$.
The classical {\bf  local Jacquet-Langlands  correspondence}   
$$JL_{\mathbb C}: {\mathcal D}_{\mathbb C}(GL(m,D))\tilde{\rightarrow} {\mathcal D}_{\mathbb C}(GL(n,F))$$
is a bijection   characterized by a  Harish-Chandra character relation on matching elliptic regular conjugacy classes.  
 Fixing an isomorphism   $\mathbb C\simeq \mathbb Q_\ell^{ac}$,
the local Jacquet-Langlands  correspondence  
  gives an $\ell$-adic  local 
Jacquet-Langlands  correspondence   $$JL_{\mathbb Q_\ell^{ac}}: {\mathcal D}_{\mathbb Q_\ell^{ac}}(GL(m,D))\tilde{\rightarrow} {\mathcal D}_{\mathbb Q_\ell^{ac}}(GL(n,F))$$
 independent of the isomorphism $\mathbb C\simeq \mathbb Q_\ell^{ac}$,
and respecting integrality. Minguez and S\'echerre \cite{MS17} proved that two integral  representations of ${\mathcal D}_{\mathbb Q_\ell^{ac}}(GL(m,D)) $ are congruent  modulo $\ell$ if and only if their  transfers to $GL(n,F)$  are congruent modulo $\ell$.  But there is no $\ell$-modular   local Jacquet-Langlands  correspondence   compatible with the $\ell$-adic  local Jacquet-Langlands  correspondence  by reduction modulo $\ell$, as for example, when $d=2$ and $q+1\equiv 0$ modulo $\ell$, the trivial representation $1_{\mathbb Q_\ell^{ac}}$ of $D^*$ corresponds to the Steinberg $\St_{\mathbb Q_\ell^{ac}}$ of $GL(2,F)$ of reduction modulo $\ell$ of length $2$ (Dat \cite{D12a}).  
 However, the   Badulescu morphism \cite{Ba07}
$$LJ_{\mathbb Q_\ell^{ac}} :{\mathcal Gr}_{\mathbb Q_\ell^{ac}}(GL(n,F))\to {\mathcal Gr}_{\mathbb Q_\ell^{ac}}(GL(m,D))$$ 
where  ${\mathcal Gr}_{R}(G)$ is the Grothendieck group  of finite length $R$-representations of $G $, gives by reduction  
an $\ell$-modular Badulescu   morphism
$$LJ_{\mathbb F_\ell^{ac}}: {\mathcal Gr}_{\mathbb F_\ell^{ac}}(GL(n,F)) \to {\mathcal Gr}_{\mathbb F_\ell^{ac}}(GL(m,D)).$$
  S\'echerre and Stevens \cite{SeSte19}  
 introduced  the  interesting  notions of  mod $\ell$ inertial supercuspidal support   and    linkage  for  irrreducible complex representations     $\pi, \pi'$ of   $G $.  
 
  a) Picking an isomorphism  $\mathbb C\simeq \mathbb Q_\ell^{ac}$ one  supposes that $\pi$ is an  $\ell$-adic representation of $G$.  The inertial class of the  cuspidal support of  $\pi$ contains an integral cuspidal representation  $\tau$. The {\bf mod $\ell$ inertial supercuspidal support } of  $\pi$ is 
the inertial  supercuspidal support of  any irreducible component of $r_\ell(\tau)$;  it depends only on the isomorphism class of $\pi$.

b)    $\pi, \pi'$  are {\bf linked} if there are  prime numbers $\ell_1, \ldots, \ell_r$  different from $p$, and     irrreducible complex representations $\pi=\pi_0, \pi_1,\ldots,\pi_r=\pi'$ such that, for each $i\in \{1,\ldots, r\}$, the representations $\pi_{i-1}, \pi_i$  have the same mod $\ell_i$ inertial supercuspidal support. 
  
When  $\pi,\pi'$ are essentially square integrable, they are linked if and only if their images  by  the local Jacquet-Langlands  correspondence  $JL_{\mathbb C}$ are linked if and only if  
 (Andrea Dotto \cite{Do21a})  they have the same semi-simple endoclass  (a  type invariant). When $G=GL(n,F)$ and $\pi,\pi'$ are cuspidal, they have the same {\bf endoclass} if and only if 
  the associated  irreducible representations of Weil group  $W_F$   by the local Langlands correspondence share an irreducible component when restricted to the wild inertia group. 
 
 \section{Bernstein   blocks}   For a commutative ring $R$, a non trivial   idempotent $e$ in   the Bernstein center $ \mathcal Z_R(G)$ decomposes  the abelian category $$\Mod_R(G) =e(\Mod_R(G)) \times (1-e)(\Mod_R(G)) $$ 
into a direct product of two abelian full subcategories. When the  idempotent  $e\in  \mathcal Z_R(G)$  is  primitive,
the subcategory $e(\Mod_R(G))$  where $e$ acts by the identity, is   indecomposable (no non trivial factors) and called  a {\bf block}.

   Bernstein and Deligne   factorized  $\Mod_{\mathbb C}(G)$ into blocks, and  described  the center of each block (a finite $\mathbb C$-algebra).  Their arguments are valid
for any algebraically closed field $R$ of characteristic $0$. We have the decomposition in blocks
 $$\Mod_R(G) =\prod_{ \Omega \in \mathcal B_R(G)} \Mod_R(G)_{\Omega}.$$ 
  over the connected components $\Omega$ of 
 the Bernstein variety   $\mathcal B_R(G)$. The {\bf Bernstein block}  $\Mod_R(G)_{\Omega}$ consists  of  the $R$-representations of $G$  all of whose irreducible  subquotients have  inertial supercuspidal support  $\Omega$. The    center of   $\Mod_R(G)_{\Omega}$ is the ring of  regular functions on the  variety $\Omega$.
The decomposition is based on the unicity of the supercuspidal support.

\bigskip  When 
 $G$ is an inner form  of $GL(n,F)$,  two complex discrete series  of $G$ in the same  block  are inertially equivalent. 
 The  Jacquet-Langlands correspondence commutes with twisting by characters,  and yields a bijection between the blocks containing discrete series.  Andrea Dotto \cite{Do21a} parametrized these  blocks  by two algebraic invariants (one is 
   the endo-class)  and  obtained a complete algebraic description of the  Jacquet-Langlands correspondence at the level of inertial classes.

 For  an algebraically closed field  $R$ of characteristic  different from $p$, the  Deligne-Bernstein decomposition   remains true (S\'echerre and Stevens \cite{SeSte16}) and Bastien Drevon and Vincent S\'echerre  \cite{DS21} described the block decomposition of the  abelian category  of finite length  $R$-representations of  $G$. Unlike the  case of all  $R$-representations of $G$, several non-isomorphic supercuspidal supports may correspond to the same block. A supercuspidal block is equivalent to the principal block of the multiplicative group of a suitable division algebra.

\bigskip   When 
 $R$ is an algebraically closed field  of characteristic banal $\ell$ for $G$,  it is expected that the Deligne-Bernstein decomposition remains true and that    the reduction modulo $\ell$ gives a bijection between the blocks of $\ell$-adic representations of $G$ and the blocks of mod $\ell$  representations of $G$.

\bigskip

 When $R= W(\mathbb F_\ell^{ac}) $  is  the Witt ring of $\mathbb F_\ell^{ac}$ and $G=GL(n,F) $, Helm \cite{Helm16a},\cite{Helm16b},\cite{Helm20} showed that   the block decomposition of $\Mod_{\mathbb F_\ell^{ac}}(G) $  lifts  to a  block decomposition of  $\Mod_{W(\mathbb F_\ell^{ac})}(G)$:
$$\Mod_{W(\mathbb F_\ell^{ac})}(G) =\prod_{ \Omega \in \mathcal B_{\mathbb F_\ell^{ac}(G)}} \Mod_{\mathbb Z_\ell^{ac}}(G)_{\Omega}.$$ 
The block $\Mod_{W(\mathbb F_\ell^{ac})}(G)_{\Omega}$ consists of the $W(\mathbb F_\ell^{ac})$-representations of $G$ such that any  irreducible subquotient  $V$ 

- has a supercuspidal support in $\Omega$ modulo isomorphism,  if $\ell V=0$

- such that the reduction modulo $\ell$ of an integral element  in the inertial class of the supercuspidal support  of $V$ is in $\Omega$  modulo isomorphism, if $\ell V=V$.

\noindent  The   center   of $\Mod_{ W(\mathbb F_\ell^{ac})}(G)_\Omega$ is a finitely generated, reduced, $\ell$-torsion free $W(\mathbb F_\ell^{ac})$-algebra  and the center of $\Mod_{W(\mathbb F_\ell^{ac})}(G)$ is naturally
   isomorphic to the ring of endomorphisms  of  the Gelfand-Graev representation of $G$ \footnote{$ \ind_U^{GL(n,F)} \psi$, where 
 $\psi $  is a generic $W(\mathbb F_\ell^{ac})$-character   of the unipotent radical $U$ of   a Borel subgroup of $GL(n,F)$} . 
  
\bigskip The blocks of $\Mod_R(G)$  have been computed in a large number of examples \footnote{with the theory of types}.  The  {\bf principal   block} of $\Mod_{R}(G)$ is the block containing  the   trivial representation of $G$. When $R=\mathbb C$, the principal block is equivalent to the category of modules over the Iwahori Hecke $\mathbb C$-algebra.  Many blocks of $\Mod_R(G)$  are equivalent to the principal block of another group $G'$.  
 When $R=\mathbb Q_\ell^{ac}, \mathbb Z_\ell^{ac}$ or $\mathbb F_\ell^{ac}$, Dat explained the known coincidences  between blocks of $\Mod_R(G)$ and predicted many more  by a functoriality principle for blocks involving dual groups \cite{Dat17}, \cite{D18a}.

{\it Example} For an algebraically closed field $R$ of characteristic different from $p$ and $G$  an inner form of $GL(n,F)$,  each   block of $\Mod_R(G)$  is equivalent to the principal  block  of a product of general linear groups \cite{SeSte16}.

 \bigskip For a commutative ring $R$ where $p$ is invertible, there is decomposition  of $\Mod_R(G)$ by the Moy-Prasad {\bf depth} (\cite{D09} Appendix A).  
 
 An $R$-representation $V$ of $G$  has {\bf depth $0$} if  $V=\sum_x V^{\tilde G_{x}}$   is the sum of its invariants $V^{\tilde G_{x}}$ by the pro-$p$ radicals of the  subgroups of $G$ fixing the vertices of the adjoint Bruhat-Tits building of $G$. 
  The possible depths   form a sequence  of non-negative rational numbers $r_0=0<r_1 <\ldots $.  The category  $\Mod_R(G)^{(r)}$ 
  of $R$-representations of 
  $G$  of depth $r$   is abelian  with an explicit  finitely generated projective generator but   is generally not a block. We have $$\Mod_{R}(G)= \prod_{n\in \mathbb N} \Mod_R(G)^{(r_n)}.$$

\bigskip  

When $p=0$ in  $R$,  the Bernstein center of  $G$ is as small as possible, equal to  the Bernstein center of the  center $Z(G)$ of $G$ (Ardakov-Schneider   \cite{ArdSchn20}  when  $R$ is a field but their proofs are  valid for a commutative ring, see also Dotto \cite{Do21a})
$$\mathcal  Z_R(Z(G)) =\varprojlim_K R[Z(G)/K ], \ \ K\subset Z(G) \ \text{open compact subgroup} .$$

 When  $E/\mathbb Q_p$ a finite extension of  ring of integers $O_E$,  the category of locally finite  (equal to the union of their subrepresentations of finite length)
  representations of $GL(2,\mathbb Q_p)$ on $O_E$-torsion modules with a central character decomposes as  a product of blocks with a noetherian  center (Paskunas and Shen-Nin Tung  \cite{PT21}).

\section{Satake isomorphism}

 The structure of the Hecke ring of any special parahoric subgroup $K$ of $G$ is understood via  the {\bf Satake transform}:
$$Sat:\mathcal H (G,K) \to \mathcal H (Z,Z^0) \ \ \ \ Sat(f)(z) =\sum_{u \in U^0\backslash U}  f(uz) \  \text{for} \ z\in Z.$$
It is an injective ring homomorphism, and as  $\mathcal H (Z,Z^0) \simeq \mathbb Z[Z/Z^0]$   is  commutative, it shows that  the Hecke ring $ \mathcal H (G,K)$ is  commutative.
A basis of the  image of  $Sat$ is  $$S_\lambda =\sum _{\lambda'\in W(\lambda)} \delta ^{1/2} (\lambda/ \lambda')e_{\lambda' } \ \ \text{for $\lambda\in Z^+/Z^0$}, $$ where $e_\lambda\in \mathcal H (Z,Z^0)$ corrresponds to $\lambda$ (Henniart-V. \cite{HV15}, \cite{HV12} Proposition 2.3). This shows that modulo isomorphism, the  commutative Hecke ring $ \mathcal H (G,K)$  does not depend on the choice of $K$.

\bigskip By scalar extension to a commutative ring $R$, the Satake transform extends  to a map $Sat:\mathcal H_R (G,K) \to  \mathcal H_R (Z,Z^0)$.
 For $R=\mathbb C$,  it is well known that  $\delta_B^{1/2} Sat$  induces an   isomorphism  $$ \mathcal H _{\mathbb C}(G,K) \simeq \mathbb C[Z/Z^0] ^{W_G}.$$  An all-important special case was singled out by Langlands, that is where $G$ is unramified and where $K$ is a hyperspecial maximal compact subgroup of $G$. Langlands interpreted the Satake isomorphism  as giving a parametrization of  the isomorphism classes of complex smooth irreducible representations of $G$ with non-zero $K$-fixed vectors, by certain semisimple conjugacy classes in a complex group $ \hat G$  ``dual'' to $G$.

 For a field $R$ of characteristic $p$,  $Sat$ induces an isomorphism (Henniart-V.\cite{HV15})\footnote{with $Z^-/Z^0$ instead of $Z^+/Z^0$ but these monoids are isomorphic}
$$ \mathcal H_R (G,K) \simeq  R[Z^+/Z^0].$$

 \bigskip Instead of focusing on  the trivial $R$-representation $1_K$ of $K$, 
 we consider   two  finitely  generated   $R$-representations $W,W'$ of $K$ and  the  Hecke $R$-bimodule  $$\mathcal H_R(G,K, W,W') \simeq \Hom_{R[G]}(\ind^G_K W, \ind^G_K W') .$$
 It  is realized as a set of compactly supported functions  $f:G\to \Hom_R(W,W')$ with a certain $K$-bi-invariance.    In the case $W=W'$, it 
  is  an algebra called an {\bf Hecke algebra with weight $W$} that we rather write  $\mathcal H_R(G,K,W)$; the Hecke algebra with  trivial weight   is the Hecke $R$-algebra $\mathcal H_R(G,K)$.
  The  Satake transform generalizes for any standard parablic subgroup $P=MN$, to an injective   map 
 $$Sat_M:\mathcal H_R(G,K, W,W') \to \mathcal H_R(M,M^0, W_{N^0},W'_{N^0}), \ \ Sat_M(f)(m) (\overline v)=\sum_{n \in N^0\backslash N} \overline { f(nm)(v) }$$
 for $m\in M, v\in W$, where $v\to \overline v$ is the quotient map $W\to W_{N^0}$ (similarly for  $W'\to W'_{N^0}$).  Another generalization considered in  (Herzig \cite{H11} when $G$ is split,  Henniart-V. \cite{HV15})   $$Sat'_M:\mathcal H_R(G,K, W,W') \to \mathcal H_R(M,M^0, W^{N^0},W'^{N^0}), \ \ Sat'_M(f)(z) ( v)=\sum_{u \in U^0\backslash U} f(uz)(v) $$
 for $v\in  W^{N^0}$,  is related to $Sat_M$  by taking  duals \cite{HV12}.  The functional approach of  $Sat_M$   (Henniart-V. \cite{HV12} Section 2) is a motivation to prefer it.

 \bigskip When $R$ is an algebraically closed field of characteristic $p$  and $W,W'$ are irreducible, the  generalized Satake transform play a role in  the modulo $p$ and $p$-adic  Langlands  correspondence. In this situation   $W_{U^0},W'_{U^0}$  have dimension $1$, the Hecke bimodule 
 $\mathcal H_R(G,K, W,W')$
  is non-zero if and  only if the $R$-characters of $Z^0$ on $W_{U^0},W'_{U^0}$  are $Z$-conjugate.  For $M=Z$,  there are explicit bases  $(S_\lambda^{W,W'})$ of the  image of  $Sat_Z$,  and 
 $(T_\lambda^{W,W'} )$ of $\mathcal H_R(G,K, W,W')$ such that 
 $$Sat_Z(T_\lambda^{W,W'})=S_\lambda^{W,W'}$$
  for $\lambda \in Z^+(W,W')/Z^0$ where $Z^+(W,W')$ is a certain union of cosets of $Z^0$ in $Z^+$ (Abe-Herzig-V. \cite{AHV22}). 
  The proof  relies on the theory the pro-$p$-Iwahori Hecke $R$-algebra. A  simple consequence   is the  ``change of weight'' \footnote{the change of weight theorem is an isomorphism between two compactly induced representations} which is an important step in the proof of the classification of admissible irreducible $R$-representations of $G$. There is  also a  change of weight the  pro-$p$-Iwahori Hecke algebra giving another proof for the change of weight for $G$  (Abe\cite{Abe17}).
 For an Hecke algebra $\mathcal H_R(G,K, W )$ with irreducible weight $W$, one gets    an explicit inverse of 
the  Satake isomorphism  (Henniart-V.\cite{HV12})  \footnote{this isomorphism for $Sat'$ is proved when $G$  is split in Herzig \cite{H11}, and in general in Henniart-V. \cite{HV15}}:
 $$Sat_Z:  \mathcal H_R(G,K, W) \ \tilde{\rightarrow}  \  \mathcal H_R( Z^+,Z^0, W_{U^0}).$$
 For $G$  quasi-split,  $ \mathcal H_R( Z^+,Z^0, W_{U^0}) \simeq R[Z^+/Z^0]$ hence  $\mathcal H_R(G,K, W )$ does not depend on the choice of 
   $(K,W)$ modulo isomorphism.  

For $G$  general,  the center of  $\mathcal H_R(G,K, W) $ contains a finitely generated  subalgebra  $\mathcal Z_T $ isomorphic to  $ R[T^+/T^0]$, and   $\mathcal H_R(G,K, W) $ is a finitely generated    $\mathcal Z_T $-module.  
One chooses   an   element $s$  in the center of $M$  which   strictly contracts $N$ by conjugation. There is a unique element $T_s\in  \mathcal H_R(M,M^0, W_{N^0})$ with support $M^0s$ such that $T_s(s)$ is the identity on  $W_{N^0}$. The generalized Satake transform  $$ Sat_M:\mathcal H_R(G,K, W) \hookrightarrow \mathcal H_R(M,M^0, W_{N^0} ) $$
is a localization at  $T_s $ \footnote{This means that  the image of $ Sat_M$ contains $T_s$ and that  its localisation at $T_s$ is $\mathcal H_R(M,M^0, W_{N^0})$.}.
The natural intertwiner 
$$I_V:\ind_K^G W\to \ind_P^G (\ind_{M^0}^M W_{N^0})$$
is injective and  its localization at $T_s$ is bijective  when $W$  satisfies a regularity assumption \footnote{meaning that the map  $\mathcal H_R(M,M^0, V_{N^0}) \otimes_{ \mathcal H_R(G,K, V) } \ind_K^G V\to \ind_P^G (\ind_{M^0}^M V_{N^0})$ is bijective, if  the kernel of  $V\to V_{N^0}$ contains $k V^{(N^{op})^0}$  for all  $k\in K \setminus P^0 (P^{op})^0$} (Herzig \cite{Her11}, Abe \cite{Abe13}, Henniart-V.\cite{HV12}).

  \bigskip  For a field  $R$ of characteristic $p$,  the  {\bf supersingularity}   of  an admissible irreducible  $R$-representation $V$ of $G$  is  defined with the Satake homomorphism  (Abe-Henniart-Herzig-V.\cite{AHHV17}). First, 
 assuming $R$-algebraically closed,  an homomorphism from the center of an Hecke algebra $\mathcal H_R(G,K, W)$ with irreducible weight  is said to be  supersingular if it 
  does  not extend  to the center of $\mathcal H_R(M, M^0, W_{N^0})$ via the Satake homomorphism  for any 
$P \neq G$. As $V$ is admissible, there exists some  irreducible  representation $W$ of $K$ such   that $\Hom_{R[G]} (\ind^G_KW, V)  \neq 0$. If   $\Hom_{R[G]} (\ind^G_KW, V)$ as a module over  the center of  $\mathcal H_R(G, K, W)$   contains an eigenvector   with a  supersingular eigenvalue, $V$ is called  { supersingular.  This does not depend on the choice of $(K,W)$. 
 For  $R$  not algebraically closed,   $V$    is  called   supersingular if 
 some  admissible  irreducible $R^{ac}$-representation $V^{ac}$  of $G$ which is $V$-isotypic as an $R$-representation,  is supersingular. This does not depend on the choice of $V^{ac}$ (Henniart-V. \cite{HV19}).

 \bigskip   
For  $G$ unramified and $K$ hyperspecial, using the geometric Satake equivalence, Xinwein Zhu \cite{Zhu21a}  identified  the Hecke ring 
$\mathcal H(G,K)$  with a ring associated to  the Vinberg monoid  of $\hat G$ and    formulated a canonical Satake isomorphism, and   proved that the commutative $\mathbb Z$-algebra $\mathcal H(G,K)$ is finitely generated. 
He extended his formulation to an Hecke  algebra $\mathcal H_{O_E}(G,K,W)$ with weight   a  finite free  $O_E$-module $W$  arising from an irreducible algebraic representation $E\otimes_{O_E}W$ of $G$, where $E/F$ is a finite extension.

  For $F $ of characteristic $0$ and $R$ a field of  characteristic $p$, Heyer  (\cite{Hey22} Theorem 4.3.2) defined  a derived Satake homomorphism related to $Sat_M$.   
  
  For $F $ of characteristic $0$,  $G$ split,  $K $  hyperspecial and   $R=\mathbb Z/p^a\mathbb Z, a\geq 1$, 
 Niccolo Ronchetti \cite{Ron19} established a Satake homomorphism for  the derived  Hecke $\mathbb Z/p^a\mathbb Z$-algebra of $(G,K)$ (a graded associative $\mathbb Z/p^a\mathbb Z$-algebra whose degree $0$ subalgebra is $\mathcal H_{\mathbb Z/p^a\mathbb Z}(G,K)$). The relation with the Heyer derived Satake homomorphism is unclear.

\section{Pro-$p$ Iwahori Hecke ring}  
 
 The isomorphism classes of the  Iwahori Hecke ring  $\mathcal H(G,J)$  and the pro-$p$ Iwahori Hecke ring    $\mathcal H(G,\tilde J)$ depend only on $G$, because  the Iwahori subgroups of $G$ are conjugate, as well as the  pro-$p$ Iwahori subgroups.
       
        They are both natural generalisations of  affine Hecke $\mathbb Z$-algebras. We will focus on the  pro-$p$ Iwahori Hecke ring which is more involved, that we will denote by $\mathcal H (G)$,  but all the results  apply to  Iwahori Hecke rings with some simplifications.

 Our motivation to study  the pro-$p$ Iwahori Hecke ring instead of the Iwahori Hecke ring  comes from the theory of  mod $p$ representations\footnote{Flicker \cite{F11} studied the pro-$p$ Iwahori Hecke complex algebra when $G$ is unramified}. Any non-zero mod $p$ representation of $G$ has a    non-zero $\tilde J$-fixed vector,   and the pro-$p$ radical of any parahoric subgroup of $G$ is contained in some $G$-conjugate of $\tilde J$.

\bigskip   For any  commutative ring $R$,    the pro-$p$ Iwahori Hecke $R$-algebra  $\mathcal H_R(G)=R\otimes_{\mathbb Z}\mathcal H(G)$  
 is a specialization of  the {\bf generic  pro-$p$ Iwahori Hecke  $R[{\bf q}_*]$-algebra}  $\mathcal H(G)({\bf q}_*,  c_* )$ of $G$, introduced by Nicolas Schmidt \cite{Schmidt09}, \cite{Schmidt17},  when $G$ is split, and by V.\cite{V16a} in general.  The $ {\bf q_*}$  are finitely many indeterminates   and   the finitely many $ c_*\in R[{\bf q}_*]$ satisfy simple conditions.  The general principle is that one proves  properties of the generic  pro-$p$ Iwahori Hecke  $R[{\bf q}_*]$-algebra by specializing all ${\bf q}_*$ to $1 $, and 
then one transfers them to  $\mathcal H_R(G)$ by specialization.   

{\it Example} The  affine Yokonuma-Hecke algebra defined by Maria Chlouveraki and Loic Poulain d'Andecy
is a generic pro-p Iwahori Hecke algebra (Chlouveraki and  S\'echerre \cite{CS16}).

  The main features \footnote{ the Iwahori Matsumoto presentation, the Bernstein basis, the Bernstein-Lusztig relations,  the description of the center, and the geometric proofs of  G\"ortz \cite{Go07}} of affine Hecke $R$-algebras generalize  to the   generic pro-$p$ Iwahori Hecke  $R$-algebra, and by specialization to $\mathcal H_R(G)$.  The  $R[{\bf q}_*]$-module $\mathcal H(G)({\bf q}_*,  c_* )$
  is  free  with an Iwahori-Matsumoto   basis of elements  satisfying  braid relations and quadratic relations,  with  ``alcove walk bases''  associated to the Weyl chambers satisfy   product formulas involving different alcove walk bases, and with  Bernstein-Lusztig relations from which one deduces  an explicit canonical $R[{\bf q}_*] $-basis  of  the center \cite{V14}. 
    
\bigskip  Finiteness properties of the pro-$p$ Iwahori ring $\mathcal H(G)$:

 (i)  The center $\mathcal Z (G)$ of  $\mathcal H(G)$  is a finitely generated $\mathbb Z$-algebra and  $\mathcal H(G)$  is a  finitely generated $\mathcal Z(G)$-module.
 
 (ii)   $\mathcal Z (G)$ contains a  canonical subring $\mathcal Z_T$ isomorphic to the affine semi-group $\mathbb Z$-algebra $\mathbb Z[T^+/ T^0]$, and 
 the $\mathcal Z_T$-modules $\mathcal Z $ and   $\mathcal H$  are  finitely generated. 

 (iii) The elements of the  Iwahori-Matsumoto basis\footnote{The Iwahori-Matsumoto basis  of $\mathcal H(G)$ is given by  the characteristic functions of the double cosets of $ G$ modulo $\tilde J$.} of   $\mathcal H(G)$ are invertible in  $\mathbb Z[1/p]\otimes_\mathbb Z \mathcal H(G)$.

(iv)  For any commutative ring $R$,  the center of $ \mathcal H_R(G) $ is  $\mathcal Z_R(G)=R\otimes_{\mathbb Z} \mathcal Z(G)$.
 
\bigskip For any field $R$, any simple $ \mathcal H_R(G)$-module is finite dimensional by (i) and (iv) \cite{H09}.

 \bigskip   Xuhua He and Radhika Ganapathy  \cite{HG21} gave  an Iwahori-Matsumoto  presentation  of the Hecke ring  ${\mathcal H} (G,J_n)$ of the $n$-th congruence subgroup $J_n$ of $J$ for any $n\in \mathbb N_{>0}$.

\bigskip For a  standard  parabolic subgroup $P=MN $,  although $\mathcal H_R(M)$ is not contained in  $\mathcal H_R(G)$,  there  is
 a {\bf parabolic induction} $$\ind_{\mathcal H(M)}^{\mathcal H(G)}= - \otimes_{\mathcal H_R(M)} X_{G,P}:\Mod  \mathcal H_R(M) \to \Mod  \mathcal H_R(G) , \ \ \ X_{G,P}=\ind_P^G ( R[\tilde J_M \backslash M])$$ 
of  right adjoint $\Hom_{\mathcal H_R(G)}(X_{G,P}, -)$  and of  left adjoint   a certain localisation (hence the left adjoint is exact,  a surprise when $p$ is not invertible in $R$ as the functor  $(-)_N$ for representations is not exact).
The parabolic induction  and its right adjoint  for the group and  for the pro-$p$ Iwahori Hecke algebra correspond to each other  via the pro-$p$ Iwahori invariant functors. 
The same holds true for the left adjoint functor  if  $R$ is a field of characteristic different from $p$, but Abe gave a counter-example  for $G=GL(2,\mathbb Q_p)$ and $R$ of characteristic $p$ (Ollivier-V.\cite{OV18}). 
The parabolic induction is  isomorphic to:    $$\ind_{ \mathcal H(P)}^{ \mathcal H(G)} =   - \otimes_{ \mathcal H(P)}  \mathcal H (G):\Mod  \mathcal H_R(M) \to \Mod  \mathcal H_R(G)  $$ where  $\mathcal H(P)=\mathbb Z[(\tilde J \cap P)\backslash G / (\tilde J \cap P)]$  is   the parabolic pro-$p$ Iwahori Hecke ring  
  of $P$ for two  ring  homomorphisms
$ \mathcal H(M) \leftarrow \mathcal H(P) \rightarrow  \mathcal H(G)$ (Claudius Heyer \cite{Hey21}).

 \bigskip 
 For an algebraically closed field   $R $ of characteristic $p$ 
and an irreducible 
 $R$-representation $W$ of a special parahoric subgroup $K$ containing $\tilde J$,  an  inverse Satake-type  isomorphism
  $$f:\mathcal H_R(Z^-,Z^0,W^{U^0}) \tilde{\rightarrow}   \mathcal H_R(G,K,W)$$
is obtained by composition of two natural algebra isomorphisms (Ollivier \cite{O15} when $G$ is split,  V.\cite{V15a} in general). The first   isomorphism is associated to a ``good'' alcove walk basis   
 $$\mathcal H_R(Z^-,Z^0,W^{U^0}) \tilde{\rightarrow}  \End_{\mathcal H_ R(G) }(W^{\tilde J} \otimes_{\mathcal H_ R(K,\tilde J)}\mathcal H_ R(G) ).$$  
The  dimension of $W^{\tilde J}$ 
is  $1$. The second  isomorphism $$
 \End_{\mathcal H_ R(G) }(W^{\tilde J} \otimes_{\mathcal H_ R(K,\tilde J)}\mathcal H_ R(G) ) \tilde{\rightarrow}  \mathcal H_R(G,K,W) $$ is associated to a natural  $H_ R(G)$-module isomorphism
  $W^{\tilde J} \otimes_{\mathcal H_ R(K,\tilde J)}\mathcal H_ R(G)  \tilde{\rightarrow}   (\ind_K^G W)^{\tilde J}$. 
 
  When $G$ is split,   $f$ is the inverse of the variant $Sat'_Z$ of the generalized Satake isomorphism (Ollivier \cite{O15}).

\section{Modules of pro-$p$ Iwahori Hecke algebras in characteristic $p$}

There is a numerical mod $p$ local Langlands correspondence for the pro-$p$ Iwahori Hecke algebra of $GL(n,F)$  (V.\cite{V05}):  the following two sets 
 
a)  the isomorphism classes of the  $n$-dimensional irreducible $\mathbb F_p^{ac}$-representations of  $\Gal(F^{ac}/F)$ with a fixed value of the determinant of the action of a Frobenius.

b) the isomorphism classes of the supersingular simple modules $ {\mathcal H}_{\mathbb F_p^{ac}}(GL(n,F))$ with a fixed action of  $p_F $  embedded diagonally. 

 \noindent have the  (finite) number of elements \footnote{equal to the number of   irreducible unitary poiynomials of degree $n$ in $k_F[X]$}.
This   was significantly improved by  Grosse-Kloenne if  
 $F \supset \mathbb Q_p$ \cite{GK16a},   \cite{GK16b}. He   constructed  an exact and fully faithful functor from the category of finite length supersingular  $ {\mathcal H}_{\mathbb F_{p^d}}(GL(n,F))$-modules to the category of $\mathbb F_q^{d}$-representations of  $\Gal(F^{ac}/F)$, if $p^d\geq q$ \footnote{$F^{sep}=F^{ac}$ as the characteristic of $F$ is $0$}.

  \bigskip We recall that  
 the pro-$p$ Iwahori Hecke ring $\mathcal H(G)$ of $G$ is a finitely generated module over a  central subring $\mathcal Z _{T}\simeq \mathbb Z[T^+/T^0]$.

 A nonzero (right) $\mathcal H_R(G)$-module $\mathcal V$ is called 
  
   {\bf ordinary} if the action on  $\mathcal V$  of  any $z\in \mathcal Z _{T}$  corresponding to a   non-invertible element  of the semi-group $T^+/ T^0$ is invertible.
 
  {\bf supersingular} if for any $v\in \mathcal V$ and    any 
 $z\in \mathcal Z _{T}$ corresponding to a non-invertible  element  of  $T^+/ T^0$,  there exists $n\in \mathbb N$ such that 
$z^n v=0$.  

\bigskip  Let $R$ be an algebraically closed field of characteristic $p$.

{\bf Classification}  
The  supersingular simple $\mathcal H_R(G)$-modules are classified (V.\cite{V15a}). 
The simple $ {\mathcal H}_R(G)$-modules  are classified in terms of the  simple  supersingular  $ {\mathcal H}_R(M)$-modules  for the Levi subgroups $M$ of the  parabolic subgroups of $G$  (Noriyuki Abe \cite{Abe19d},  algebraically closed is not necessary (Henniart-V.\cite{HV19}):  
  
    For a  standard parabolic subgroup $P=MN$ of $G$ and a simple supersingular $ \mathcal H_R(M) $-module $\mathcal W$, there is a notion of extension $e_{P'}(\mathcal W)$ of $\mathcal W$ to $ \mathcal H_R(M') $ for  a parabolic subgroup $P'=M'N'$ of $G$ containing $P$.  There is a maximal  $P'$  with this property, denoted by $P(\mathcal W)$. For a parabolic subgroup $Q$ with $P\subset Q\subset P(\mathcal W)$, there is    a generalized Steinberg    $ \mathcal H_R(M(\mathcal W)) $-module
 $$ st_Q^{P(\mathcal W)} (\mathcal W)=\ind_{\mathcal H(Q)} ^{\mathcal H(G)} (e_Q(\mathcal W))/\sum_{Q\subsetneq Q'\subset Q(\mathcal W)}ind_{\mathcal H(Q')} ^{\mathcal H(G)} (e_{Q'}(\mathcal W)).$$
  The triple $(P,\mathcal W, Q)$ is called standard.  
 The $\mathcal H_R(G)$-module  $$I_{\mathcal H(G)}(P,\mathcal  W,Q)= \ind_{\mathcal H(P(\mathcal W))} ^{\mathcal H(G)}  (st_Q^{P(\mathcal W)} (\mathcal W))$$  is simple, and any simple  $\mathcal H_R(G)$-module is isomorphic to  $I_{\mathcal H(G)}(P,W,Q)$ for some standard triple $(P,W,Q)$ unique modulo $G$-conjugation.  It is 
 ordinary if and only if $P=B$.

 \bigskip {\bf Extensions} The extensions between simple $ \mathcal H_R(G)$-modules $$\Ext^i _{ \mathcal H_R(G)}   (I_{\mathcal H_R(G)}(P_1,\mathcal W_1,Q), I_{\mathcal H(G)}(P_2,\mathcal  W_2,Q_2)),  \ \ \ i\geq0,$$
are either  $0$, or 
extensions  between  supersingular simple modules of a  specialization of a  generic pro-$p$ Iwahori Hecke algebra which is not of  a pro-$p$ Iwahori Hecke $R$-algebra (Abe \cite{Abe22}).  
In more details, considering the central characters, 
 the extensions are $0$ if $P_1\neq P_2$. When $P=P_1=P_2$, following the construction of the simple modules, 
 $$\Ext^i _{ \mathcal H_R(G)} (I_{\mathcal H(G)}(P,\mathcal W_1,Q), I_{\mathcal H(G)}(P,\mathcal  W_2,Q_2))  \simeq \Ext^i_{\mathcal H_R(M')}( 
 st^{P'}_{{\mathcal Q}_1'}(\mathcal W_1) ,  st^{P'}_{\mathcal Q_2'}(\mathcal W_2)) \ $$
for some   $P', Q'_1, Q'_2$,
$$ \Ext^i_{\mathcal H_R(G)}( 
 st^{G}_{{\mathcal Q}_1}(\mathcal W_1) ,  st^{G}_{\mathcal Q_2}(\mathcal W_2))\simeq  \Ext^{i-r}_{\mathcal H_ R(G)}( 
 e_G(\mathcal W_1) ,  e_G(\mathcal W_2)),$$
 for some explicit $r\in \mathbb N_{\geq 0}$, and using results of Ollivier-Schneider \cite{OS14},  
 $$  \Ext^{i}_{\mathcal H_ R(G)}( 
 e_G(\mathcal W_1) ,  e_G(\mathcal W_2)) \simeq \Ext^i_{\mathcal H_ R(M)/I}(\mathcal W_1, \mathcal W_2) $$
 for some ideal $I$ of $\mathcal H_ R(M)$ acting on $\mathcal W_1, \mathcal W_2$ by $0$.   Abe   computed explicitely 
  $\Ext^1 $ for  two supersingular simple  $\mathcal H_ R(M)/I$-modules.

\bigskip  
 When    $G = GL(2,F)$, C\'edric P\'epin and Tobias Schmidt  proved:
 
 (i)  The $2$-dimensional supersingular simple $\mathcal H_{\mathbb F_p^{ac}}(G)$-modules   can be realized through the equivariant cohomology of the flag variety of the dual group $\hat G$ over  $\mathbb F_p^{ac}$ \cite{PS19}. 
 
 (ii) There is 
    a version in families of  the Breuil's semisimple mod $p$ Langlands correspondence for $GL_2(\mathbb Q_p)$   \cite{PS20}. 
    
    (iii) There is  a 
Kazhdan-Lusztig theory for the  generic pro-$p$ Iwahori  Hecke $\mathbb Z[{\bf q}]$-algebra  of  $ G$, where the role of   $\hat G$ is taken by the Vinberg monoid $V_{\hat G}$ and its flag variety; the monoid comes with a fibration $V_{\hat G}\to \mathbb A^1$ and the dual parametrisation of $\mathcal \mathcal H_{\mathbb F_p^{ac}(}G)$-modules is achieved by working  
over the $0$-fiber. They  introduce  the  generic pro-$p$ antispherical module and  the  generic pro-$p$  Satake homomorphism  for a generic spherical Hecke $\mathbb Z[{\bf q}]$-algebra \cite{PS21}.

\section{Representations  in characteristic  $p$} 
 In this section,   $R$  is a  field of   characteristic $p$.  The admissible irreducible $R$-representations of $G$ are classified  in terms of the supersingular admissible irreducible $R$-representations of the Levi subgroups of $G$ (Abe-Henniart-Herzig-V. \cite{AHHV17}   for $R$ algebraically closed, Henniart-V. \cite{HV19}  
 for $R$ not  algebraically closed).
 
  \bigskip {\bf Classification} The representation $\ind_P^GW$ parabolically induced from an  irreducible admissible supersingular  $R$-representation $W$ of a  Levi subgroup $M$ of a parabolic subgroup $P$ of $G$,  has multiplicity $1$ 
and irreducible subquotients  
 $$I_G(P,W,Q)= \ind_{P(W)}^G (e(W)\otimes \St_Q^{P(W)})$$  for the parabolic subgroups $Q$ of $G$ containing $P$ and contained in the maximal parabolic subgroup  $P(W)$ where 
  the inflation of  $W$   to $P$  extends to a representation $e(W)$, and 
   $$\St_Q^{P(W)}=(\ind^{P(W)}_Q 1_{Q})/  \sum_{Q\subsetneq Q' \subset P(W)} \ind^{P(W)}_{Q'} 1_{Q'}.$$ 
   Any irreducible admissible $R$-representation $V$ of $G$ is isomorphic to $I_G(P,W,Q)$ for a unique   triple $(P,W,Q)$  modulo $G$-conjugation.

   \bigskip  A similar classification holds true for the irreducible admissible genuine mod $p$ representations of the metaplectic double cover of $Sp_{2n}(F)$ (Koziol and Laura Peskin \cite{KP18}).

   \bigskip There is  a complete description of $\ind_P^GW$ for any irreducible admissible $R$-representation $W$ of $M$  \cite{HV19}.
As a corollary, one obtains generic irreducibility and for any admissible irreducible $R$-representation $V$ of $G$,

 \centerline{ $V$ supersingular  $\Leftrightarrow $  $V$   cuspidal $\Leftrightarrow$  $V $  supercuspidal.  }

  The supersingularity of  $V$  can also be seen on  
the  pro-$p$ Iwahori invariants  (Ollivier-V.\cite{OV18}   for $R$ algebraically closed,  but algebraically closed is not necessary Henniart-V.\cite{HV19}):

$V$ is supersingular $\Leftrightarrow $ $V^{\tilde J}$ is supersingular $\Leftrightarrow$ some  non-zero subquotient of $V^{\tilde  J}$ is supersingular.

   The   $\tilde J$-invariant of  $I_G(P,W,Q)$ of $G$ is computed and depends only on 
the pro-$p$ Iwahori invariant $W^{\tilde J_M}$ (Abe-Henniart-V.\cite{AHV18}, \cite{AHV19}).

\bigskip  When $P=B$, the irreducible subquotients of $\ind_B^G W$ are called  {\bf ordinary}.  An admissible $R$-representation of $G$ with  ordinary irreducible subquotients is called ordinary. 

The $\tilde J$-invariant functor   induces an equivalence between the category of finite length ordinary $R$-representations of $G$ generated by their $\tilde J$-invariant vectors   and the category  of the finite length ordinary $ \mathcal H_R(G)$-modules, assuming $R$ algebraically closed (Abe  \cite{Abe19b}).

\bigskip  
When  $F$ has  characteristic $0$, the higher duals $(S^i( I_G(P,W,Q)))_{i\geq 0}$  are computed in terms of 
  $(S^i  (W))_{i\geq 0}$  in  a few cases (Kohlhaase \cite{Kohl17}). 
  
The extensions between $R$-representations $\ind_P^G W$ of $G$ parabolically induced from  supersingular absolutely irreducible $R$-representations  $W$ of   Levi subgroups,  are computed in many cases   when $G$ is split and $R$ finite (Hauseux \cite{Hau16},\cite{Hau17},\cite{Hau18a},\cite{Hau19}).

 \bigskip The  supersingular  admissible irreducible   $R$-representations $V$ of $G$ are not understood, this remains an open crucial question since two decades 
and  a stumbling block for the search of a mod $p$ or $p$-adic local Langlands correspondence if $G\neq GL(2,\mathbb Q_p)$. 
  The  classification of simple supersingular $\mathcal H_R(G)$-modules   does not help because we  do not have enough information on the   pro-$p$ Iwahori invariant functor.

\bigskip 
  Breuil \cite{Br03} relying on the work of Barthel-Livne   classified the  supersingular  admissible irreducible mod $p$ representations of  $GL(2,\mathbb Q_p)$.  This  was the starting point of 
  the mod $p$ local Langlands correspondence  for $GL(2,\mathbb Q_p)$.    There are two main novel features in the mod $p$ local Langlands correspondence. It involves reducible    representations of $GL_2(\mathbb Q_p)$,   and it  extends to an exact functor from finite length  representations of $GL_2(\mathbb Q_p)$  to  finite length representations of $\Gal(\mathbb Q_p^{ac}/\mathbb Q_p)$. 
  
\bigskip There are scarce results on supersingular admissible irreducible  mod $p$ representations of  $G\neq GL(2,\mathbb Q_p)$.  Supersingular admissible irreducible  mod $p$ representations  are  classified only for  some groups close to $GL(2,\mathbb Q_p)$:  for $SL(2,\mathbb Q_p)$  (Ramla Abdellatif \cite{Abde14}, Chuangxun Cheng \cite{Che13}), and for the unramified unitary group $U(1,1)(\mathbb Q_p^2/\mathbb Q_p)$ in two variables (Koziol \cite{Koz16a}).

 When $F\neq \mathbb Q_p$, there can be  many more  supersingular admissible irreducible  mod $p$ representations of $GL(2,F)$
   than $2$-dimensional irreducible representations of $\Gal(F^{sep}/F)$  (Breuil-Paskunas \cite{BP12}); they cannot be described as quotients of a compact induction by a finite number of equations (Hu \cite{Hu12}   if $F\supset \mathbb F_p((T))$, Schraen \cite{Schraen15}   if $F/ \mathbb Q_p$ is quadratic, Wu \cite{Wu21}  in general  if $F\supsetneq \mathbb Q_p$).  
   
   When $R$ is a field of characteristic $p$ and $F\supset  \mathbb Q_p$, Herzig-Koziol-V.\cite{HKV20} proved 
  that any $G$ admits a supersingular admissible irreducible  $R$-representation, using a local method of Paskunas  \cite{Pask04} if the semisimple rank $r_G$ of $G$ is $1$, and a global method  if  $r_G>1$.  The existence is not known if  $F\supset \mathbb F_p((T))$. 
  
 \bigskip    
There have   been  recent advances which strongly suggest that the study of mod $p$  representations of $G$ is best done on the derived level.  When $R$ is a field of characteristic $p$, Schneider \cite{Schn15} introduced  the unbounded derived category $D_R(G)$ of $R$-representations of $G$. When  $\tilde J$ is  torsion free (this forces $F$ to be of characteristic $0$),  $D_R(G)$ is equivalent  to the  derived category of differential graded modules over a differential graded version $\mathcal H_R(G)^\bullet$  of the pro-$p$ Iwahori Hecke $R$-algebra of $G$, by the derived  $\tilde J$-invariant functor.

The parabolic induction $\ind_P^G:\Mod_R(M)\to \Mod_R(G)$ being exact extends to   an exact derived parabolic induction $R\ind_P^G: D_R(M)\to D_R(G)$ between the unbounded derived categories. The total derived functor of $R_P^G$   is right adjoint to  $R\ind_P^G$. The category $D_R(G)$ has arbitrary  small direct products  and $R\ind_P^G$ commutes with arbitrary  small direct products (Heyer \cite{Hey22}), hence $R\ind_P^G$ has a left adjoint\footnote{by Brown representability; Heyer \cite{Hey22} gave  another proof}. When $\tilde J$ is  torsion free, the functor $R\ind_P^G$  corresponds  to the  derived parabolic induction functor on the dg Hecke algebra side, via the  derived  $\tilde J$-invariant functor  (Sarah  Scherotzke and Schneider \cite{SS22}).  

The Kohlhaase duality functors are related to the derived duality functor $\RHom (-,R)$   (Schneider and Claus Sorensen \cite{SchnSo22}). 
 
  The structure of the cohomology algebra $\Ext^\bullet_{\Mod_R(G)} (R[\tilde J \backslash G], R[\tilde J \backslash G])$  is simpler than   of $\mathcal H_R(G)^\bullet$;  there is an explicit description when $G=SL(2,\mathbb Q_p), p\geq 5$ (Ollivier and Schneider \cite{OS19},  \cite{OS21}).

     \section{Local Langlands correspondences  for $GL(n,F)$}

  The classical local Langlands  correspondence for $GL(n,F)$ is   a  bijection between 
  the  isomorphism classes of irreducible   complex representations of $GL(n,F)$ and the isomorphism classes of $n$-dimensional  Weil-Deligne  complex representations,
  given by local class field theory when $n=1$, and characterized by  the  requirement  that the $L$ and $\epsilon $ \footnote{for a fixed non trivial $\mathbb C$-character of $F$} factors attached to   corresponding pairs of complex representations  coincide (Henniart \cite{H02}).    
 A $n$-dimensional  Weil-Deligne  complex representation  is a pair $(\sigma , N)$ where $\sigma $ is a $n$-dimensional Frobenius semi-simple  complex representation of the Weil group $W_F$  and $N\in \End_{\mathbb C}\sigma$ is a nilpotent   endomorphism satisfying $\sigma(w)N\sigma(w)^{-1}=q^{|w|} N$ for all $w\in W_F$ \footnote{$|w| $  is the  power of $q$ to which $w$ raises the elements of the residue field $k_F$}.
   The supercuspidal irreducible   $\mathbb C$-representations of $GL(n,F)$ correspond to  the $n$-dimensional irreducible    $\mathbb C$-representations of  $W_F$  \footnote{$N=0$}.

  A twist of the correspondence by an unramified complex character of $GL(n,F)$  is  compatible  with all the automorphisms of $\mathbb C$. For a prime $r$, an isomorphism $\mathbb C\simeq \mathbb Q_r^{ac}$  transfers  the  twisted classical local Langlands correspondence  to a local Langlands correspondence  for $\mathbb Q_r^{ac}$-representations of $GL(n,F)$.
For $r=\ell\neq p$, the nilpotent part is related to the action of the tame inertia group on an  $\ell$-adic representation of  $W_F$.
By  reduction modulo $\ell$  of the   $\ell$-adic 
local Langlands correspondence composed with the Zelevinski involution on $\ell$-adic representations of $GL(n,F)$, one obtains a {\bf $\ell$-modular local correspondence}.  The   $\ell$-modular local   correspondence,  is a parametrization for $\ell$-modular irreducible representations of $GL(n,F)$ by $n$-dimensional  Weil-Deligne  $\ell$-modular representations, defined as above with $\mathbb F_\ell^{ac}$ instead of $\mathbb C$.   The supercuspidal irreducible   $\ell$-modular representations of $GL(n,F)$ correspond to  the $n$-dimensional irreducible   $\ell$-modular representations of  $W_F$. But the nilpotent part $N$  of the Weil-Deligne $\ell$-modular representation has   no obvious Galois interpretation.

\bigskip  Dat   \cite{D12a}, \cite{D12b}, \cite{D12c}   obtained  a geometric realization  of the $\ell$-adic local Zelevinski correspondence and  of the   $\ell$-modular  local  correspondence  on the unipotent\footnote{$=$ in the principal block $=$ subquotients of some $\Ind_B^G(\chi)$ for $\chi$ an unramified character of a Borel subgroup $B$, this is not the definition of Lusztig} irreducible $\mathbb F_\ell^{ac}$-representations of $GL(n,F)$   when the order of $q$ in $\mathbb F_\ell^*$ is at least $n$ \cite{D12d}\footnote{the regular case}, and  when $q\equiv 1\mod \ell$ and $\ell >n$  \cite{D14}\footnote{the limit case}.

\bigskip  Kurinczuk and Matringe  \cite{KM17},\cite{KM19a},\cite{KM19b},\cite{KM20a}, extended to   $\ell$-modular   representations  the Rankin-Selberg local constants of Jacquet, Piatetski-Shapiro and Shalika of  pairs of complex generic representations of linear groups, and the Artin-Deligne local constants of  pairs of  complex Weil-Deligne representations. These local constants are preserved by the classical local Langlands correspondence, but not by 
the $\ell$-modular local  correspondence.  Enlarging the space of  $\ell$-modular Weil-Deligne representations to representations with non necessarily nilpotent operators, they suggested a {\bf modified   $\ell$-modular local   corrrespondence}  
compatible with the formation of local constants and  characterized by a list of natural properties.  When  $R$  is a noetherian $W(\overline F_\ell ^{ac})$-algebra, using the Rankin-Selberg functional equations, Matringe and Moss \cite{MM22} proved that  a  $R$-representation of  $GL(n,F)$ of Whittaker type admits  a Kirillov model.
  
 \bigskip When the characteristic of $F$ is $0$, Breuil and Schneider \cite{BS07} motivated  by a $p$-adic extension of the local Langlands  correspondence,  suggested a {\bf modified    local Langlands correspondence}  where  the complex representations of $GL(n,F)$ are no more irreducible. 
 The Langlands  quotient theorem realizes an  irreducible  $\mathbb C$-representation $V$ of $GL(n,F)$ as a quotient of a certain parabolically induced representation  $\ind_P^GW$.  In the modified version,  $V$ is replaced by  by a twist of  $\ind_P^GW$ by an unramified character of $GL(n,F)$.
 
\bigskip  When the characteristic of $F$ is $0$,  Emerton and Helm \cite{EH14} motivated by a local Langlands correspondence in families and  by global contexts,   introduced the {\bf generic  $\ell$-adic local Langlands  correspondence} 
  which  has useful applications to the cohomology of  Shimura varieties.   For any finite extension $E/\mathbb Q_\ell$, it is a map  $\rho \mapsto \pi(\rho)$  from $n$-dimensional continuous $E$-representations  of the Galois group
 $\Gal(F^{ac}/F)$ to finite length $E$-representations  of $GL(n,F)$ with  an absolutely irreducible  generic socle and no other generic  irreducible subquotients\footnote{It is a slight modification of the Breuil and Schneider correspondence transferred to $\ell$-adic representations; the socle of $V$ is the maximal semi-simple subrepresentation of $V$}. 
 Each  
$ \pi(\rho)$ contains a $GL(n,F)$-stable $O_E$-lattice $\pi(\rho)^o$ of reduction modulo $\ell$ having an absolutely irreducible socle and no other generic  subquotients, unique modulo homotethy.

The {\bf generic mod $\ell$ local Langlands  correspondence} (Emerton-Helm \cite{EH14})   is 
a  modified mod $\ell$   local  correspondence compatible  with  the generic $\ell$-adic local Langlands correspondence by reduction modulo $\ell$. Irreducible representations of $GL(n,F)$ are no longer irrreducible,  Weil-Deligne representations are now  Galois representations,  and the Zelevinski involution does not intervene. For a  finite extension $R/\mathbb F_\ell$, it is  the unique  map $\overline \rho \mapsto \overline \pi(\overline \rho)$  from $n$-dimensional  $R$-representations  of  $\Gal(F^{ac}/F)$ to   finite length $R$-representations  of $GL(n,F)$ such that 
 
 1) $\overline \pi(\overline \rho)$ has an absolutely irreducible  generic socle and no other generic  irreducible subquotients,  

 2) For all finite extensions $E/\mathbb Q_\ell$  of ring of integers $O_E$ and residue field $k_E$ containing $R$, and $\rho: \Gal(F^{ac}/F)\to GL(n,O_E)$  a continuous representation lifting $\overline \rho \otimes_R k_E$  the  the reduction modulo $\ell$ of  $\pi(\rho)^o$ \footnote{$\rho$ identifies with a representation $\Gal(F^{ac}/F)\to GL(n,E)$} embeds in $\overline \pi(\overline \rho)\otimes_R k_E$. 
 
 3)   $ \overline \pi(\overline \rho)$ is minimal with respect to the above two conditions.  
 
 \noindent For $GL(2,F)$, the correspondence   is fairly concrete and  explicit  when $\ell \neq 2$ (Helm \cite{Helm13}).

 Emerton and Helm introduced also   a notion of  a {\bf local Langlands correspondence in families} \footnote{For an example of a local p-adic Langlands correspondence  in families for $GL(2,\mathbb Q_p)$, see Ildar Gaisin and Joaquin Rodrigues Jacinto  \cite{GR17}}. For  any  suitable complete local noetherean algebra $R$ with finite residue field $k$, it is the unique  map $\rho \mapsto \pi(\rho)$ from the continuous representations $\rho:\Gal(F^{ac}/F)\to GL(n,R)$ to the  
 admissible $R$-representations  of $GL(n,F)$ that interpolates the generic local Langlands correspondence over the points of $\Spec R$ and satisfies certain technical hypotheses.

  The existence of the map (Helm and Moss  \cite{HM18})  amounts to showing that whenever there is a congruence between  two $\ell$-adic representations of $\Gal(F^{ac}/F)$, there is a corresponding congruence on the other side of the $\ell$-adic generic local Langlands correspondence. The key idea of the proof is the introduction of    the  Bernstein center $\mathcal Z$  of $\Mod_{ \mathbb Z_\ell^{ac} }(GL(n,F))$ (Helm \cite{Helm16a}, \cite{Helm16b}, \cite{Helm20}). It encodes deep information about  congruences between $\mathbb Z_\ell^{ac}$-representations of $GL(n,F)$. For instance, if two  integral  irreducible $\mathbb Q_\ell^{ac}$-representations of $GL(n,F)$ become isomorphic modulo $\ell$,  then $\mathcal Z$  acts on these representations  by scalars congruent modulo $\ell$.

 \bigskip  For any $R$  and  $G$ quasi-split,  Dat, Helm, Kurinczuk and Moss \cite{DHKM20}  studied  the scheme of Langlands parameters  of $G$ with  coefficient   the smallest possible ring $R=\mathbb Z[1/p]$  for    a local Langlands correspondence in families. In particular, this allows to study chain of congruences of Langlands parameters modulo several different primes.  
The  blocks in  the category of   $\mathcal Z_{\mathbb Z^{ac} [1/p]}$-representations of $G$   of   depth $0$ are in natural bijection with the connected components of the  space of tamely ramified Langlands parameters for $G$ over $\mathbb Z^{ac} [1/p]$; the category  is indecomposable  if $G$  semi-simple and simply connected, or  unramified (Dat and Thomas Lanard \cite{DL22}).
In a work in progress, Dat, Helm, Kurinczuk and Moss  extend the Emerton-Helm-Moss    local Langlands correspondence in families  to a conjecture which asserts  the existence of   isomorphisms between  

a) the centre of  $\Mod_{\mathbb Z[q^{-1/2}]}(G)$,  

b) the ring of  functions on the moduli stack of Langlands parameters\footnote{Constructions of moduli spaces of Langlands parameters have been also proposed by  Fargues and Scholze (\cite{FS21} over $\mathbb Z_\ell, \ell \neq p$ using the condensed mathematics of Clausen-Scholze, and by  Xinwen Zhu \cite{Zhu21b}). The   local Langlands correspondence  is now conjectured to  exist at a categorical level (Denis Gaitsgory \cite{Ga16}).} for $G$ over $\mathbb Z[q^{-1/2}]$, 

c) the descent to $\mathbb Z[q^{-1/2}]$ of the endomorphisms of a Gelfand-Graev representation of $G$. 

\noindent They prove this conjecture when $G$ is any classical $p$-adic group after inverting an integer. The conjecture   follows from  a  Fargues-Scholze conjecture  (\cite{FS21} I.10.2) \footnote{private communication of Dat}. 

     \bigskip The $p$-adic local Langlands correspondence for $GL (n,F)$, $F$ of characteristic $0$,   is a hypothetical correspondence between continuous unitary $E$-Banach space representations of  $GL(n,F)$   and $n$-dimensional continuous $E$-representations  of $\Gal(F^{ac}/F)$, for any finite extension   $E/\mathbb Q_p$,  given by local class field theory  when $n=1$,  as  all local Langlands corespondences for $GL (n,F)$. Ana Caraiani, Matthew Emerton, Toby Gee, David Geraghty, Paskunas and Sug Woo Shin \cite{CEGGPS16} constructed a  candidate  when $p $ does not divide $2n$ using global methods.     For  $F=\mathbb Q_p$ and $n=2$, it coincides with  the $p$-adic local correspondence envisionned by Breuil twenty years ago, constructed by Pierre Colmez \cite{Co10},  and analyzed  by  (Paskunas \cite{Pask13}, Colmez, Gabriel Dospinescu, Paskunas \cite{CDP14}).

  For $n\geq 2$ and  $D_n$  the central division algebra  
  over $F$ of invariant $1/n$, Scholze\cite{Scho18} constructed a candidate for a $p$-adic  and mod $p$ Jacquet-Langlands correspondence from $GL(n,F) $ to $D_n^*$   in a purely geometric way, using the cohomology of the infinite-level Lubin-Tate space.   The  mod $p$ Jacquet-Langlands correspondence is a canonical  map from the  admissible mod $p$ representations of $GL(n,F)$
to the  admissible mod $p$ representations of $ D_n^*$ having  a continous action of $\Gal(F^{ac}/F)$.
For $F=\mathbb Q_p$ and $n=2$,  it is studied in  (Dospinescu-Paskunas-Schraen \cite{DPS22}).

\section{Gelfand-Kirillov Dimension}

Let $R$ be a field and $V$  an  irreducible admissible $R$-representation of $G$. 
 For any decreasing sequence $(K_i)_{i\geq 1}$ of open compact subgroups of $G$  with limit the trivial group,  the  dimensions $\dim_RV^{K_i}$ for $i\geq 1$ are finite and form an increasing sequence.
 If   $V$ is finite dimensional, $\dim_RV^{K_i}=\dim_R V$ when $i$ is large enough.
 Generally, the dimension of $V$ is infinite and   $(\dim_RV^{K_i})_{i\geq 1}$ tends to infinity,   but how  ?

 \bigskip  For  $F$ of characteristic $0$, one can choose   an   $O_F$-lattice $L$ of the Lie algebra $\mathcal G$ of $G$ on which  the exponential map is defined,  and consider the decreasing sequence $(K_i=\exp(p_F^{2i} L))_{ i\geq 1}$.  
       When $R=\mathbb C$, the  Harish-Chandra local character expansion of $V$ implies  that  $\dim_R V^{K_i}$ becomes eventually polynomial\footnote{Harish-Chandra, Notes by Stephen DeBacker and Paul J. Sally, Admissible invariant distributions on reductive $p$-adic group, University Lecture Series Vol.16,1999,97p.}
 $$\dim_R V^{K_i}=P_{L,V}(q^i),  \ \ P_{L,V} (X)\in \mathbb Q[X]  \ \ \ \ \text{
 for $i$ large enough. } $$    The degree  $d_{V}$  of the polynomial $P_{L,V}[X]$   does not depend on the choice of $L$.  It is half  the dimension of a unipotent conjugacy class in $G$,
     $$0\leq d_{V}\leq \dim_F U,$$
      and  is $0$ if and only if  $V$ is finite dimensional.  
 The  growth of $\dim_R V^{K_i}$ is measured by $q^{d_V}$.
 
\bigskip   Wth no restriction on $F$, for $G=GL(n,F)$,  $K_i=1+p_F^i M(n,O_F)$ for $i\geq 1$, and   $R$   of characteristic different from $p$,   the dimensions $\dim_RV^{K_i} $ satisfy the above properties  \footnote{article in preparation}.  
  For $R=\mathbb C$, this follows from the 
   Roger Howe  local character expansion. For  $R$  of characteristic $\ell$,  the key of the proof is that 
  any
cuspidal  irreducible $ \ell$-modular representation  of $GL(n,F)$   lifts to an  irreducible cuspidal $ \ell$-adic representation.

For $GL(2,F)$, we have $\dim_F  U=1$, and $V$ is infinite dimensional if and only if   $d_V=1$. For $GL(n,F)$, we have  $\dim_F  U=(n^2-n)/2$ and $V$ is generic (in particular if $V$ is cuspidal)  if and only if  $d_V=\dim_F  U$.
 
For $R$  a  finite field of characteristic $p$,
 $G=GL(2, \mathbb Q_p)$,  and $V$ absolutely irreducible,  the dimensions $\dim_RV^{K_i} $  for $i\geq 1$ computed by  Stefano Morra \cite{Mo13} satisfy the above properties. 

\bigskip For  $F $ of characteristic $0$, $R$  a  finite field of characteristic $p$,  $K$ a uniformly powerful open pro-$p$ subgroup of $G$,   $K_i$ the closed subgroup of $K$ generated by   $\{k^{p^i}, k\in K\}$  for $i\geq 1$,  and $V$ an admissible $R$-representation of $G$, there is  a positive integer $\delta_{V}$ not depending on the choice of $K$  and  positive real numbers $a\leq b$, such that (Calegari-Emerton\cite{CE11}, Emerton-Paskunas \cite{EP20}, Dospinescu-Paskunas-Schraen \cite{DPS22}):
 $$a p^{i \delta_{V}}  \leq \dim_R V^{K_i}\leq b p^{i \delta_{V}} . $$
    The integer $\delta_V$  which is a sort of Iwasawa dimension of the dual of $V$, is  called the Gelfand-Kirillov dimension of $V$. When  $F/\mathbb Q_p$ is unramified, the admissible  $R$-representations $V$ of $GL_2(F) $ studied by Breuil-Herzig-Hu-Morra-Schraen\cite{BHHMS21b}   in   mod $p$ cohomology  satisfy  $\delta_V=[F:\mathbb Q_p] $.  If $V$ is    isomorphic to $I_G(P,W,Q)$ we have\footnote{article in preparation} 
  $$ \delta_{ I_G(P,W,Q)}=\delta_W + \dim _{\mathbb Q_p} N_Q$$ 
where $ N_Q$ is the unipotent radical of the parabolic subgroup  $Q$ of $G$. 

\bigskip These partial results indicate that  perhaps a notion of Gelfand-Kirillov dimension  of $V$  exists  for any $F, G ,R$.

  \end{document}